\documentclass[12pt,twoside,american,english,british]{article}
\usepackage[T1]{fontenc}
\usepackage[utf8]{inputenc}
\usepackage[a4paper]{geometry}
\usepackage{color}
\usepackage{textcomp}
\usepackage{enumitem}
\usepackage{amsmath}
\usepackage{amsthm}
\usepackage{amssymb}
\usepackage{setspace}
\usepackage{esint}

\makeatletter
\numberwithin{equation}{section}
\numberwithin{figure}{section}
\theoremstyle{plain}
\newtheorem{thm}{\protect\theoremname}[section]
\theoremstyle{definition}
\newtheorem{defn}[thm]{\protect\definitionname}
\theoremstyle{plain}
\newtheorem{lem}[thm]{\protect\lemmaname}
\theoremstyle{plain}
\newtheorem{prop}[thm]{\protect\propositionname}
\theoremstyle{definition}
\newtheorem{example}[thm]{\protect\examplename}
\newcommand{\lyxaddress}[1]{
	\par {\raggedright #1
	\vspace{1.4em}
	\noindent\par}
}

\usepackage[T1]{fontenc}
\usepackage[utf8]{inputenc}
\usepackage[a4paper]{geometry}
\usepackage{color}
\usepackage{dsfont}
\usepackage{amsmath}
\usepackage{amsthm}
\usepackage{amssymb}
\usepackage{setspace}
\usepackage{esint}

\makeatletter
\numberwithin{equation}{section}
\numberwithin{figure}{section}
\theoremstyle{plain}

\@ifundefined{date}{}{\date{}}
\usepackage{babel}
\usepackage{latexsym}
\usepackage{amsthm}
\usepackage{esint}
\usepackage{eucal}
\usepackage{epstopdf}
\usepackage{graphicx}
\usepackage{bigints}
\theoremstyle{plain}

\newtheoremstyle{boldremark}
    {\dimexpr\topsep/2\relax} 
    {\dimexpr\topsep/2\relax} 
    {}          
    {}          
    {\bfseries} 
    {.}         
    {.5em}      
    {}          

\theoremstyle{boldremark}
\newtheorem{brem} [thm] {Remark} 

\usepackage{fancyhdr}
\usepackage{blindtext}
\usepackage{titlesec}

\titleformat{name=\chapter}[display]
{\normalfont\Huge\itshape}
{\titlerule[1pt]\vspace{-33pt}\filleft%
  \parbox[t]{6em}{%
    \raggedleft%
    \rule{\linewidth}{0.5ex}\newline%
    \chaptertitlename\ \thechapter%
  }%
}
{4pc}
{\normalfont\upshape\bfseries\Huge}
\allowdisplaybreaks

\makeatother

\usepackage{babel}
\addto\captionsbritish{\renewcommand{\definitionname}{Definition}}
\addto\captionsbritish{\renewcommand{\lemmaname}{Lemma}}
\addto\captionsbritish{\renewcommand{\theoremname}{Theorem}}
\addto\captionsenglish{\renewcommand{\definitionname}{Definition}}
\addto\captionsenglish{\renewcommand{\lemmaname}{Lemma}}
\addto\captionsenglish{\renewcommand{\theoremname}{Theorem}}
\providecommand{\definitionname}{Definition}
\providecommand{\lemmaname}{Lemma}
\providecommand{\theoremname}{Theorem}

\usepackage[colorlinks,pdfpagelabels,pdfstartview = FitH,bookmarksopen
= true,bookmarksnumbered = true,linkcolor = red,plainpages =
false,hypertexnames = false,citecolor = blue,pagebackref=false,urlcolor=blue]{hyperref}

\makeatother

\usepackage{babel}
\addto\captionsamerican{\renewcommand{\definitionname}{Definition}}
\addto\captionsamerican{\renewcommand{\examplename}{Example}}
\addto\captionsamerican{\renewcommand{\lemmaname}{Lemma}}
\addto\captionsamerican{\renewcommand{\propositionname}{Proposition}}
\addto\captionsamerican{\renewcommand{\theoremname}{Theorem}}
\addto\captionsbritish{\renewcommand{\definitionname}{Definition}}
\addto\captionsbritish{\renewcommand{\examplename}{Example}}
\addto\captionsbritish{\renewcommand{\lemmaname}{Lemma}}
\addto\captionsbritish{\renewcommand{\propositionname}{Proposition}}
\addto\captionsbritish{\renewcommand{\theoremname}{Theorem}}
\addto\captionsenglish{\renewcommand{\definitionname}{Definition}}
\addto\captionsenglish{\renewcommand{\examplename}{Example}}
\addto\captionsenglish{\renewcommand{\lemmaname}{Lemma}}
\addto\captionsenglish{\renewcommand{\propositionname}{Proposition}}
\addto\captionsenglish{\renewcommand{\theoremname}{Theorem}}
\providecommand{\definitionname}{Definition}
\providecommand{\examplename}{Example}
\providecommand{\lemmaname}{Lemma}
\providecommand{\propositionname}{Proposition}
\providecommand{\theoremname}{Theorem}

\begin{document}
\title{\textbf{Regularity of vectorial minimizers for non-uniformly elliptic
anisotropic integrals}}
\author{Pasquale Ambrosio, Giovanni Cupini, Elvira Mascolo}
\date{February 11, 2026 }

\maketitle

\begin{abstract}
\begin{singlespace}
\noindent We establish the local boundedness of the local minimizers
$u:\Omega\rightarrow\mathbb{R}^{m}$ of non-uniformly elliptic integrals
of the form $\int_{\Omega}f(x,Dv)\,dx$, where $\Omega$ is a bounded
open subset of $\mathbb{R}^{n}$ ($n\geq2)$ and the integrand satisfies
anisotropic growth conditions of the type 
\[
\sum_{i=1}^{n}\lambda_{i}(x)|\xi_{i}|^{p_{i}}\le f(x,\xi)\le\mu(x)\left\{ 1+|\xi|^{q}\right\} 
\]
for some exponents $q\geq p_{i}>1$ and with non-negative functions
$\lambda_{i},\mu$ fulfilling suitable summability assumptions. The
main novelties here are the degenerate and anisotropic behaviour of
the integrand and the fact that we also address the case of vectorial
minimizers ($m>1$). Our proof is based on the celebrated Moser iteration
technique and employs an embedding result for anisotropic Sobolev
spaces.\vspace{0.2cm}
\end{singlespace}
\end{abstract}

\selectlanguage{english}%
\begin{singlespace}

\section{Introduction}
\end{singlespace}

\selectlanguage{british}%
\begin{singlespace}
\noindent $\hspace*{1em}$In this paper we are interested in the regularity
of local minimizers $u:\Omega\rightarrow\mathbb{R}^{m}$,\foreignlanguage{english}{
$u\in W^{1,1}(\Omega;\mathbb{R}^{m})$ with $m\geq1$,} of non-uniformly
elliptic functionals of the form\foreignlanguage{english}{
\begin{equation}
\mathcal{F}(v)=\int_{\Omega}f(x,Dv)\,dx,\label{eq:funzionale}
\end{equation}
where $\Omega$ is a bounded open subset of $\mathbb{R}^{n}$, $n\geq2$.
We assume that the energy density $f=f(x,\xi)$, $x\in\Omega$, $\xi\in\mathbb{R}^{m\times n}$,
is a Carathéodory function, convex and of class $C^{1}$ with respect
to $\xi$ and satisfying the following \textit{degenerate }and\textit{
anisotropic }behaviour: for some exponents $p_{i}$, $i\in\{1,\ldots,n\}$,
and $q$ with $1<p_{i}\le q$ and for some measurable functions $\lambda_{i},\mu:\Omega\rightarrow[0,\infty)$,
$i\in\{1,\ldots,n\}$,
\begin{equation}
\sum_{i=1}^{n}\lambda_{i}(x)|\xi_{i}|^{p_{i}}\le f(x,\xi)\le\mu(x)\left\{ 1+|\xi|^{q}\right\} \quad\,\,\,\mathrm{{\color{red}{\normalcolor for\,\,a.e.\,\,\mathit{x}\in\Omega\,\,and\,\,for\,\,every\,\,\mathit{\xi\in\mathbb{R}^{m\times n}}}},}\label{eq:crescita}
\end{equation}
where
\begin{equation}
\lambda_{i}^{-1}\in L_{loc}^{r_{i}}(\Omega),\,\,\,{\color{red}{\normalcolor i\in\{1,...,n\}}},\;\;\;\;\;\;{\color{red}{\normalcolor \mu\in L_{loc}^{s}(\Omega)}}\label{eq:integrability}
\end{equation}
for some ${\color{red}{\normalcolor r_{i}\in[1,\infty]}}$ and $s\in(1,\infty]$.
Throughout this paper $\xi_{i}$, $i\in\{1,\ldots,n\}$, denotes the
$i$-th column of the $m\times n$-matrix $\xi=(\xi_{i}^{\alpha})$,
$i\in\{1,...,n\}$, $\alpha\in\{1,...,m\},$ i.e.
\[
\xi=(\xi_{1},\xi_{2},\ldots,\xi_{n})=\left(\begin{array}{cccc}
\xi_{1}^{1} & \xi_{2}^{1} & \cdots & \xi_{n}^{1}\\
\xi_{1}^{2} & \xi_{2}^{2} & \cdots & \xi_{n}^{2}\\
\vdots & \vdots & \ddots & \vdots\\
\xi_{1}^{m} & \xi_{2}^{m} & \cdots & \xi_{n}^{m}
\end{array}\right).
\]
In particular, if $\xi=Dv$, then $\xi_{i}=(v_{x_{i}}^{1},...,v_{x_{i}}^{m})^{T}$.
In addition, we assume that there exists a function $g:\Omega\times[0,\infty)^{n}\to[0,\infty)$
such that
\begin{equation}
f(x,\xi)\,=\,g(x,|\xi_{1}\vert,\ldots,|\xi_{i}\vert,\ldots,|\xi_{n}\vert)\label{eq:struttura-1}
\end{equation}
for a.e. $x\in\Omega$ and every $\xi\in\mathbb{R}^{m\times n}$.}\\
$\hspace*{1em}$In t\foreignlanguage{english}{he case $p:=\min\,\{p_{i}\}<q$,
with $\lambda_{i}\geq1$ for every $i$ and $\mu\in L^{\infty}(\Omega)$,
the functional (\ref{eq:funzionale}) belongs to the class of variational
problems with }$p,q$-growth \foreignlanguage{english}{conditions,
introduced by Marcellini \cite{Mar0,Mar1,Mar2,Mar3} and since then
widely investigated. As for the case of anisotropic }$p_{i},q$-growth\foreignlanguage{english}{,
we recall e.g. Boccardo-Marcellini-Sbordone \cite{BocMarSbo}, Stroffolini
\cite{Stroff}, Fusco-Sbordone \cite{FuSbo} and Marcellini \cite{Mar0.1,Mar2.1,Mar3bis}.}\\
$\hspace*{1em}$\foreignlanguage{english}{Anisotropic elliptic equations
have been considered under many different aspects, for instance with
respect to the maximum principle and the multiplicity of solutions:
see e.g. Pucci, R\v{a}dulescu\textit{ et al.} \cite{BouPuRad}, \cite{FoMuPu}
and \cite{MiPuRa}.}\\
$\hspace*{1em}$\foreignlanguage{english}{Actually, the research on
problems satisfying }$p,q$-growth \foreignlanguage{english}{conditions
is so intense that it is impossible to give an exhaustive and comprehensive
list of references; for an overview on the subject and a detailed
bibliography, see Mingione \cite{Ming}, }Marcellini \cite{Mar4}
and Mingione-R\v{a}dulesc\foreignlanguage{english}{u \cite{MingRad}.}\\
$\hspace*{1em}$\foreignlanguage{english}{In the vector-valued case,
as suggested by well-known counterexamples by De Giorgi \cite{De Giorgi},
Giusti-Miranda \cite{GiuMir}, }Ne\v{c}\foreignlanguage{english}{as
\cite{Necas} and} \v{S}ver\'{a}k-Yan\textit{ }\foreignlanguage{english}{\cite{Sverak-Yan},
the structural assumption $f(x,Dv)=F(x,\vert Dv\vert)$ on the integrand
is generally required for everywhere regularity. We point out that
the condition $f(x,Dv)=F(x,\vert Dv\vert)$ is more specific than
(\ref{eq:struttura-1})}.\\
$\hspace*{1em}$\foreignlanguage{english}{For the $L^{\infty}$ regularity
in the vectorial framework, see for example \cite{CMM2012,cupmarmas0,Dall'Aglio-Mas,LeoMas},
where the authors established the local boundedness of solutions for
some classes of quasilinear systems, which - in the variational context
- may correspond to integrals as in (\ref{eq:funzionale}). In particular,
in \cite{cupmarmas0} local $L^{\infty}$-estimates were obtained
for the local minimizers of (\ref{eq:funzionale}), by assuming (\ref{eq:struttura-1})
and (\ref{eq:crescita}) }with positive constants in place of the
functions $\lambda_{i}(x)$ and $\mu(x)$\foreignlanguage{english}{,
i.e. without imposing growth conditions on the integrand that are
degenerate in $x$.}\\
$\hspace*{1em}$\foreignlanguage{english}{In the context of non-uniform
ellipticity, starting from the celebrated paper by Trudinger \cite{Tru},
in \cite{CMM2018} the authors proved the local boundedness for local
minimizers of integrals of the form (\ref{eq:funzionale}) under}
$p,q$-growth c\foreignlanguage{english}{onditions of the type 
\[
\lambda(x)\,|\xi|^{p}\le f(x,\xi)\le\mu(x)\left\{ 1+\vert\xi|^{q}\right\} ,
\]
for some exponents $q\geq p>1$ and with non-negative functions $\lambda,\mu$
satisfying suitable summability conditions. We also refer to \cite{BiCuMa},
where the local boundedness is established for scalar-valued quasi-minimizers
of non-uniformly elliptic integrals of the form\pagebreak
\[
\mathcal{G}(v)=\int_{\Omega}\tilde{f}(x,v,Dv)\,dx\,,
\]
where $\tilde{f}:\Omega\times\mathbb{R}\times\mathbb{R}^{n}\rightarrow\mathbb{R}$
is a Carathéodory map satisfying a non-uniform growth condition of
the type
\[
\lambda(x)\,|\xi|^{p}\le\tilde{f}(x,u,\xi)\le\mu(x)\left\{ |\xi|^{p}+\vert u\vert^{q}\right\} +a(x)\,,
\]
for $q\geq p>1$ and non-negative functions $a,\lambda,\mu$ fulfilling
appropriate summability assumptions.}\\
$\hspace*{1em}$\foreignlanguage{english}{In this paper we carry on
the research conducted in the previous articles, extending the (local)
$L^{\infty}$-regularity results to the case of vectorial local minimizers
($m>1$) of non-uniformly elliptic integrals under }the anisotropic
growth conditions (\ref{eq:crescita}), which \foreignlanguage{english}{may
possibly be degenerate with respect to the $x$-variable. Precise
assumptions and statements are given in Section \ref{assumptions},
where we also impose some restrictions to the growth exponents $\{p_{i}\}$
and $q$ in (\ref{eq:crescita}) and to the integrability exponents
$\{r_{i}\}$ and $s$ in (\ref{eq:integrability}). These restrictions
are the natural generalizations of the various bounds that must be
satisfied to ensure the regularity of local minimizers; for a knowledge
of some of these bounds, we limit ourselves to pointing out the works
\cite{BiCuMa,BocMarSbo,cupmarmas1,cupmarmas0,CMM2017,CMM2018,CMM2023,CuMaMaPas,FePaPo}
and the references therein, since there is a vast literature on this
topic. In this regard, see also Remark \ref{comparison} below, where
we compare the main result of this paper (Theorem \ref{thm:theo1})
with those obtained in \cite{cupmarmas1,cupmarmas0,CMM2017} and \cite{CMM2018}.}\\
$\hspace*{1em}$\foreignlanguage{english}{The main novelty of our
Theorem \ref{thm:theo1} is the degenerate and anisotropic behaviour
of the integrand (\ref{eq:crescita}) in the \textit{vector case}
$m>1$. In this respect, we wish to mention the very recent work \cite{GAO}:
there, Feng, Gao and Zhang establish the local boundedness of vectorial
local minimizers for a class of integral functionals with rank-one
convex integrands and specific structural conditions. In particular,
their result is applicable to integrals of the type
\begin{equation}
\bigintssss_{\Omega}\left\{ \sum_{\alpha=1}^{m}\lambda(x)\,\vert Dv^{\alpha}\vert^{p}\,+\,\mu(x)\,\vert Dv\vert^{r}\right\} dx\,,\label{eq:Gao}
\end{equation}
with suitable $\lambda(x),\mu(x)>0$ and $p,r>1$. However, the result
in \cite{GAO} is not comparable to ours, since the functionals of
the form (\ref{eq:Gao}) do not satisfy conditions (\ref{eq:crescita})
and (\ref{eq:struttura-1}), not even if $n=m\geq2$, $\lambda_{1}=\cdots=\lambda_{n}=\lambda$
and $p_{1}=\cdots=p_{n}=p$.}\\
\foreignlanguage{english}{$\hspace*{1em}$Before describing the structure
of this paper, we also want to }point out the recent article \cite{FePaPo},
where\foreignlanguage{english}{ Feo, Passarelli di Napoli and Posteraro}
prove the local boundedness of \textit{scalar} minimizers of non-uniformly
elliptic integrals of the form (\ref{eq:funzionale}), assuming that
the left inequality in (\ref{eq:crescita}) holds with monomial weights
of the type
\[
\lambda_{i}(x)=\vert x_{i}\vert^{\alpha_{i}\,p_{i}}\,\,\,\,\,\,\,\,\,\mathrm{for\,\,some\,\,}\alpha_{i}\in[0,1)\,\,\,\,\,\,\,\,\,\mathrm{for\,\,every\,\,}i\in\{1,\ldots,n\}.
\]
The proof methods used in \cite{FePaPo} and \cite{GAO} adapt the
renowned De Giorgi iteration technique. Moreover, the approach employed
in \cite{FePaPo} relies essentially on an anisotropic Sobolev inequality
with respect to the weights $\vert x_{i}\vert^{\alpha_{i}\,p_{i}}$.\\
\foreignlanguage{english}{$\hspace*{1em}$}It is worth noting that,
to deal with the anisotropic behaviour of the integrand, we base our
estimates on an embedding result for anisotropic Sobolev spaces due
to Troisi \cite{Troisi}. Other key ingredients in the proof of our
result are the derivation of the Euler’s equation for the functional
(\ref{eq:funzionale}) and a suitable Moser iteration procedure.\foreignlanguage{english}{}\\
\foreignlanguage{english}{$\hspace*{1em}$}The paper is organized
as follows. In the next section we give the complete statement of
the main regularity result. Section \ref{sec:prelim} is devoted to
the preliminaries: after a list of some classic notations and some
essential lemmas, we recall the definition and properties of the anisotropic
Sobolev spaces that will be needed to prove our result. In Section
\ref{s:eulero} we prove that an Euler's equation holds true. This
is a main step in the proof of Theorem \ref{thm:theo1}, which is
given in Section \ref{sec:Proof}. Finally, in Section \ref{sec:example}
we present interesting examples of applicability of our main result.\vspace{1cm}

\end{singlespace}
\selectlanguage{english}%
\begin{singlespace}

\section{Assumptions and statement of the main result\label{assumptions}}
\end{singlespace}

\selectlanguage{british}%
\begin{singlespace}
\noindent $\hspace*{1em}$\foreignlanguage{english}{Let us define
the integral functional
\begin{equation}
\mathcal{F}(v):=\int_{\Omega}f(x,Dv(x))\,dx,\label{functional}
\end{equation}
where $\Omega$ is a bounded open subset of $\mathbb{R}^{n}$, $n\geq2$,
and $v\in W^{1,1}(\Omega;\mathbb{R}^{m})$, $m\in\mathbb{N}$.}\\
$\hspace*{1em}$\foreignlanguage{english}{In the sequel, we denote
by $\mathbb{R}_{+}$ the set $[0,\infty)$. Moreover, we assume that
$f:\Omega\times\mathbb{R}^{m\times n}\rightarrow\mathbb{R}_{+}$ is
a Carathéodory function satisfying the following conditions:\medskip{}
}
\end{singlespace}
\selectlanguage{english}%
\begin{itemize}
\item[\foreignlanguage{english}{\textbf{(A1)}}] $f=f(x,\xi)$, $x\in\Omega$, $\xi\in\mathbb{R}^{m\times n}$, is
convex and of class $C^{1}$ with respect to $\xi$;
\begin{singlespace}
\item[\foreignlanguage{english}{\textbf{(A2)}}] there exists a function $g:\Omega\times(\mathbb{R}_{+})^{n}\to\mathbb{R}_{+}$
such that
\begin{equation}
f(x,\xi)\,=\,g(x,|\xi_{1}\vert,\ldots,|\xi_{i}\vert,\ldots,|\xi_{n}\vert)\quad\mathrm{for\,\,a.e.}\,\,x\in\Omega\,\,\mathrm{and}\mathrm{\,\,every}\,\,\mathit{\xi\in\mathbb{R}^{m\times n}};\label{eq:struttura}
\end{equation}
\item[\foreignlanguage{english}{\textbf{(A3)}}] there exists $\tau\geq1$ such that 
\[
f(x,t\xi)\,\le\,t^{\tau}f(x,\xi)
\]
for every $t>1$, for a.e. $x\in\Omega$ and every\textcolor{red}{{}
}$\xi\in\mathbb{R}^{m\times n}$;
\item[\foreignlanguage{english}{\textbf{(A4)}}] there exist some exponents $p_{i}$, $i\in\{1,\ldots,n\}$, and $q$
with $1<p_{i}\le q$ and some measurable functions $\lambda_{i},\mu:\Omega\rightarrow\mathbb{R}_{+}$,
$i\in\{1,\ldots,n\}$, such that 
\begin{equation}
\sum_{i=1}^{n}\lambda_{i}(x)|\xi_{i}|^{p_{i}}\le f(x,\xi)\le\mu(x)\left\{ 1+|\xi|^{q}\right\} \quad\mathrm{{\color{red}{\normalcolor for\,\,a.e.\,\,\mathit{x}\in\Omega\,\,and\,\,every\,\,\mathit{\xi\in\mathbb{R}^{m\times n}}}},}\label{growth}
\end{equation}
where 
\[
\lambda_{i}^{-1}\in L_{loc}^{r_{i}}(\Omega),\,\,\,{\color{red}{\normalcolor i\in\{1,...,n\}}},\;\;\;\;\;\;{\color{red}{\normalcolor \mu\in L_{loc}^{s}(\Omega)}}
\]
for some ${\color{red}{\normalcolor r_{i}\in[1,\infty]}}$ and $s\in(1,\infty]$.\vskip.3cm
\end{singlespace}
\end{itemize}
\selectlanguage{british}%
\begin{singlespace}
\noindent $\hspace*{1em}$For every \foreignlanguage{english}{$i\in\{1,\ldots,n\}$,}
we denote
\begin{equation}
\sigma_{i}=\begin{cases}
\begin{array}{cc}
\frac{p_{i}r_{i}}{r_{i}+1} & \,\,\,\mathrm{if}\,\,\,r_{i}\in[1,\infty),\\
p_{i} & \mathrm{if}\,\,\,r_{i}=\infty,\,\,\,\,\,\,
\end{array}\end{cases}\label{eq:sigma_i}
\end{equation}
and additionally require that $\sigma_{i}\geq1$, i.e. we require
that $r_{i}\geq\frac{1}{p_{i}-1}$ if $1<p_{i}<2$.

\noindent $\hspace*{1em}$\foreignlanguage{english}{We set
\[
W^{1,\mathcal{F}}(\Omega;\mathbb{R}^{m}):=\{v\in W^{1,1}(\Omega;\mathbb{R}^{m})\,:\,\mathcal{F}(v)<+\infty\}
\]
and define a local minimizer of (\ref{functional}) as follows.}
\end{singlespace}
\selectlanguage{english}%
\begin{defn}
\begin{singlespace}
\noindent A function $u$ is a local minimizer of (\ref{functional})
if $u\in W^{1,\mathcal{F}}(\Omega;\mathbb{R}^{m})$ and 
\[
\mathcal{F}(u)\leq\mathcal{F}(u+\varphi)
\]
for all $\varphi\in W^{1,\mathcal{F}}(\Omega;\mathbb{R}^{m})$ with
${\color{red}{\normalcolor \mathrm{supp}\,\varphi\Subset\Omega}}$.
\end{singlespace}
\end{defn}

\selectlanguage{british}%
\begin{singlespace}
\noindent $\hspace*{1em}$\foreignlanguage{english}{Given a real number
$\ell\ge1$, $\ell^{*}$ is its Sobolev exponent, i.e.
\begin{equation}
\ell^{\ast}:=\begin{cases}
\begin{array}{cc}
\frac{n\ell}{n-\ell} & \mathrm{if}\,\,\,\ell<n,\\
\mathrm{any\,\,value\,\,in}\,\,(\ell,\infty) & \mathrm{if}\,\,\,\ell\geq n,
\end{array}\end{cases}\label{eq:SobConj}
\end{equation}
and $\ell'$ is the conjugate exponent of $\ell$, i.e. $\ell'=\frac{\ell}{\ell-1}$
if $\ell>1$, while $\ell'=\infty$ if $\ell=1$. We set 
\[
\mathbf{p}:=(p_{1},\ldots,p_{n}),\,\,\,\,\,\,\,\,\,\,\mathbf{pr}:=(p_{1}r_{1},\ldots,p_{n}r_{n}),\,\,\,\,\,\,\,\,\,\,\mathbf{\sigma}:=(\sigma_{1},\ldots,\sigma_{n})
\]
and let $\overline{\mathbf{p}}$ be the harmonic average of $\{p_{i}\}$,
i.e.
\begin{equation}
\frac{1}{\overline{\mathbf{p}}}\,=\,\frac{1}{n}\,\sum_{i=1}^{n}\frac{1}{p_{i}}\,.\label{eq:harmonic}
\end{equation}
Similarly, we will denote by $\overline{\mathbf{pr}}$ and $\overline{{\bf \sigma}}$
the harmonic average of $\{p_{i}r_{i}\}$ and $\{\sigma_{i}\}$, respectively.
As usual, $\frac{1}{\infty}$ has to be read as $0$. Therefore, in
the particular case $r_{1}=\cdots=r_{n}=\infty$, we have $\frac{1}{\overline{\mathbf{pr}}}=0$.
Likewise, $\frac{1}{qs}=0$ if $s=\infty$.}\\
$\hspace*{1em}$\foreignlanguage{english}{As we anticipated earlier,
to prove the local boundedness of the local minimizers of (\ref{functional}),
we need some restrictions on the exponents $\{p_{i}\},q,\{r_{i}\}$
and $s.$ More precisely, our main result reads as follows.}
\end{singlespace}
\selectlanguage{english}%
\begin{thm}
\begin{singlespace}
\noindent \label{thm:theo1}Let us assume that $(\mathbf{A1})$, $(\mathbf{A2})$,
$(\mathbf{A3})$ and $(\mathbf{A4})$ hold under the summability conditions
\[
\lambda_{i}^{-1}\in L_{loc}^{r_{i}}(\Omega),\,\,\,{\color{red}{\normalcolor i\in\{1,...,n\}}},\;\;\;\;\;\;{\color{red}{\normalcolor \mu\in L_{loc}^{s}(\Omega)}},
\]
\smallskip{}
for some ${\color{red}{\normalcolor r_{i}\in[1,\infty]}}$ and $s\in(1,\infty]$
such that $\max\,\{\sigma_{i}\}<\overline{\sigma}^{*}$ and
\begin{equation}
\frac{1}{\overline{\mathbf{pr}}}\,+\,\frac{1}{qs}\,+\,\frac{1}{\overline{\mathbf{p}}}\,-\,\frac{1}{q}\,<\,\frac{1}{n}\,.\label{eq:restriction}
\end{equation}
\foreignlanguage{british}{For every }$i\in\{1,\ldots,n\}$,\foreignlanguage{british}{
}if\foreignlanguage{british}{ $1<p_{i}<2$} we also require\foreignlanguage{british}{
$r_{i}\geq\frac{1}{p_{i}-1}$}. Then, every local minimizer $u$ of
$(\ref{functional})$ is locally bounded. Moreover, for every ball
$B_{R_{0}}(x_{0})\Subset\Omega$ with $R_{0}\in(0,1]$,\foreignlanguage{british}{}\\
\foreignlanguage{british}{}\\
\hspace*{1em}\foreignlanguage{british}{$\mathrm{(}1\mathrm{)}\,\,$there}\hspace{0.4em}\foreignlanguage{british}{exists}\hspace{0.4em}\foreignlanguage{british}{a}\hspace{0.4em}\foreignlanguage{british}{constant}\hspace{0.4em}\foreignlanguage{british}{$c_{1}$}\hspace{0.4em}\foreignlanguage{british}{$>$}\hspace{0.4em}\foreignlanguage{british}{$1$,}\hspace{0.4em}\foreignlanguage{british}{depending}\hspace{0.4em}\foreignlanguage{british}{on}\hspace{0.4em}\foreignlanguage{british}{$m$,}\hspace{0.4em}\foreignlanguage{british}{$n$,}\hspace{0.4em}\foreignlanguage{british}{$\mathbf{p}$,}\hspace{0.4em}\foreignlanguage{british}{$q$,}\hspace{0.4em}\foreignlanguage{british}{$\mathbf{r}$,}\hspace{0.4em}\foreignlanguage{british}{$s$,}\hspace{0.4em}\foreignlanguage{british}{$\tau$}\hspace{0.4em}\foreignlanguage{british}{if}\hspace{0.4em}\foreignlanguage{british}{$\overline{\sigma}$}\hspace{0.4em}\foreignlanguage{british}{$<$}\hspace{0.4em}\foreignlanguage{british}{$n$,}\hspace{0.4em}\foreignlanguage{british}{and}\hspace{0.4em}\foreignlanguage{british}{also}\\
\hspace*{1em}\foreignlanguage{british}{$\,\,\,\,\,\,\,\,\,$on $R_{0}$
if $\overline{\sigma}\geq n$, such that for every $R\in(0,R_{0})$
we have\begin{align}\label{eq:u-bound1}
\Vert u\Vert_{L^{\infty}(B_{R/2}(x_{0}))}\,\leq\,\frac{c_{1}}{R^{\vartheta_{3}}}\left[1+\Vert\mu\Vert_{L^{s}(B_{R})}^{\frac{1}{\overline{\mathbf{p}}}}\prod_{i=1}^{n}\Vert\lambda_{i}^{-1}\Vert_{L^{r_{i}}(B_{R})}^{\frac{1}{np_{i}}}\right]^{\vartheta_{1}}\big\Vert\vert u\vert+1\big\Vert_{L^{qs'}(B_{R})}^{\vartheta_{2}}\,,
\end{align} }\hspace*{1em}\foreignlanguage{british}{$\mathrm{(}2\mathrm{)\,}\,$there
exists a constant $c_{2}>0$, depending on $m$, $n$, $\mathbf{p}$,
$q$, $\mathbf{r}$, $s$, $\tau$ and $R_{0}$, such that for}\\
\hspace*{1em}\foreignlanguage{british}{$\,\,\,\,\,\,\,\,\,$every
$R\in(0,R_{0})$ we have\begin{align}\label{eq:u-bound2}
\Vert u-u_{R}\Vert_{L^{\infty}(B_{R/(2\sqrt{n})}(x_{0}))}\,&\leq\,\frac{c_{2}}{R^{\vartheta_{3}}}\left[1+\Vert\mu\Vert_{L^{s}(B_{R/\sqrt{n}})}^{\frac{1}{\overline{\mathbf{p}}}}\prod_{i=1}^{n}\Vert\lambda_{i}^{-1}\Vert_{L^{r_{i}}(B_{R/\sqrt{n}})}^{\frac{1}{np_{i}}}\right]^{\vartheta_{1}}\nonumber\\
&\,\,\,\,\,\,\,\,\,\,\,\,\,\,\,\,\cdot\left[1+\left(1+\int_{B_{R}(x_{0})}f(x,Du)\,dx\right)^{\frac{1}{p}}\,\sum_{i=1}^{n}\Vert\lambda_{i}^{-1}\Vert_{L^{r_{i}}(B_{R})}^{\frac{1}{p_{i}}}\right]^{\vartheta_{2}},
\end{align}where $u_{R}:=\frac{1}{\vert B_{R}(x_{0})\vert}\int_{B_{R}(x_{0})}u\,dx$,
$p:=\underset{1\,\leq\,i\,\leq\,n}{\min}\{p_{i}\}$ and 
\[
\vartheta_{1}:=\,\frac{\overline{\sigma}^{*}}{\overline{\sigma}^{*}-qs'}\,,\,\,\,\,\,\,\,\,\,\,\vartheta_{2}:=\,\frac{q\,(\overline{\sigma}^{*}-\overline{\mathbf{p}}s')}{\overline{\mathbf{p}}\,(\overline{\sigma}^{*}-qs')}\,,\,\,\,\,\,\,\,\,\,\,\vartheta_{3}:=\,\frac{\overline{\sigma}^{*}\,[q^{2}s'+n\,(q-\overline{\mathbf{p}})]}{\overline{\mathbf{p}}qs'\,(\overline{\sigma}^{*}-qs')}\,.
\]
}
\end{singlespace}
\end{thm}

\selectlanguage{british}%
\noindent \begin{brem}\label{thm:qs'}If $n>\overline{\sigma}$,
then condition (\ref{eq:restriction}) is equivalent to 
\begin{equation}
qs'<\overline{\sigma}^{*}.\label{eq:restriction2}
\end{equation}
In the case $n\leq\overline{\sigma}$, we can still suppose that $qs'<\overline{\sigma}^{*}$,
provided we choose a sufficiently large value for $\overline{\sigma}^{*}$.
For our purposes, in the sequel we will always assume that $qs'<\overline{\sigma}^{*}$.

\selectlanguage{english}%
\noindent $\hspace*{1em}$Also note that if $p=q$ (i.e. \foreignlanguage{british}{$p_{1}=\cdots=p_{n}=q$),
then $\vartheta_{2}=1$ and $\vartheta_{3}=\vartheta_{1}$.\end{brem}}

\selectlanguage{british}%
\noindent \begin{brem}[\textbf{Comparison with other results}]\label{comparison}In
the previous paper \cite{CMM2018}, the authors assume that there
exist measurable functions $\lambda,\mu:\Omega\rightarrow\mathbb{R}_{+}$
such that
\begin{equation}
\begin{cases}
\begin{array}{c}
\lambda(x)\,|\xi|^{p}\le f(x,\xi)\le\mu(x)\left\{ 1+\vert\xi|^{q}\right\} \\
\lambda^{-1}\in L_{loc}^{r}(\Omega),\;\;\;{\color{red}{\normalcolor \mu\in L_{loc}^{s}(\Omega)}}\,\,\,\,\,\,\,\,\,\,\,\,\,\,\,\,\,\,
\end{array}\end{cases}\label{eq:p,q-growth}
\end{equation}
for $1<p\leq q$, for some exponents\foreignlanguage{english}{ $r\in[1,\infty]$,
$s\in(1,\infty]$} and for every $\xi\in\mathbb{R}^{m\times n}$.
Therefore, the left-hand side of (\ref{growth}) is an anisotropic
version of the left inequality in (\ref{eq:p,q-growth}). Furthermore,
condition (\ref{eq:restriction}) in Theorem \ref{thm:theo1} is the
evident counterpart of condition $(2.6)$ in \cite[Theorem 2.1]{CMM2018},
which reads as follows:
\[
\frac{1}{pr}\,+\,\frac{1}{qs}\,+\,\frac{1}{p}\,-\,\frac{1}{q}\,<\,\frac{1}{n}\,.
\]
In this regard, also observe that (\ref{eq:restriction2}) is a generalization
of the restriction $q<\overline{\mathbf{p}}^{*}$ adopted, for example,
in \cite{cupmarmas1,cupmarmas0,CMM2017}, where the authors assume
that $\mu$ is a positive constant and $\lambda_{i}\equiv1$ for any
$i\in\{1,\ldots,n\}$. Moreover, the requirement \foreignlanguage{english}{$\max\,\{\sigma_{i}\}<\overline{\sigma}^{*}$
}in Theorem \ref{thm:theo1} is the natural generalization of the
condition \foreignlanguage{english}{$\max\,\{p_{i}\}<\overline{\mathbf{p}}^{*}$
}in \cite[Theorem 2.3]{CMM2017}. In fact, if $r_{1}=\cdots=r_{n}=\infty$,
then $\sigma_{i}=p_{i}$ for every $i\in\{1,\ldots,n\}$, and therefore
$\overline{\sigma}=\overline{\mathbf{p}}$.\end{brem}
\selectlanguage{english}%
\begin{singlespace}

\section{Notations and preliminaries\label{sec:prelim}}
\end{singlespace}

\begin{singlespace}
\noindent $\hspace*{1em}$In this paper we shall denote by $C$ or
$c$ a general positive constant that may vary on different occasions.
Relevant dependencies on parameters and special constants will be
suitably emphasized using parentheses or subscripts. \foreignlanguage{british}{The
norm we use on $\mathbb{R}^{k}$, }\foreignlanguage{american}{$k\in\mathbb{N}$}\foreignlanguage{british}{,
will be the standard Euclidean one and it will be denoted by $\left|\,\cdot\,\right|$.
In particular, for the vectors $\mathbf{v},\mathbf{w}\in\mathbb{R}^{k}$,
we write $\langle\mathbf{v},\mathbf{w}\rangle$ for the usual inner
product and $\left|\mathbf{v}\right|:=\langle\mathbf{v},\mathbf{v}\rangle^{\frac{1}{2}}$
for the corresponding Euclidean norm.}\\
$\hspace*{1em}$In what follows, $B_{r}(x_{0})=\left\{ x\in\mathbb{R}^{n}:\left|x-x_{0}\right|<r\right\} $
will denote the $n$-dimensional open ball with radius $r>0$ and
center $x_{0}\in\mathbb{R}^{n}$. We shall sometimes omit the dependence
on the center when all balls occurring in a proof are concentric.
Unless otherwise stated, different balls in the same context will
have the same center.

\selectlanguage{british}%
\noindent $\hspace*{1em}$\foreignlanguage{american}{If $E\subseteq\mathbb{R}^{k}$
is a Lebesgue-measurable set, then we will denote by $\vert E\vert$
its $k$-dimensional Lebesgue measure. When $0<\vert E\vert<\infty$,
the mean value of a function $v\in L^{1}(E)$ is defined by 
\[
\fint_{E}v(x)\,dx\,:=\,\frac{1}{\vert E\vert}\int_{E}v(x)\,dx\,.
\]
}$\hspace*{1em}$Now we gather some results that will be needed later
on. Let us start with a lemma on the properties of the functions $f$
and $g$ considered in Section \ref{assumptions}.\medskip{}

\end{singlespace}
\begin{lem}
\begin{singlespace}
\noindent \label{lem:g function}Assume \foreignlanguage{english}{$(\mathbf{A1})$
and} \foreignlanguage{english}{$(\mathbf{A2})$. Then, for a.e. $x\in\Omega$
and every }$i\in\{1,\ldots,n\}$,\foreignlanguage{english}{\medskip{}
}

\selectlanguage{english}%
\noindent \textup{(i)}$\,\,\,g(x,t_{1},\ldots,t_{n})$ \foreignlanguage{british}{is
of class $C^{1}$ with respect to $(t_{1},...,t_{n})\in(\mathbb{R}_{+})^{n}$};\medskip{}

\noindent \textup{(ii)}$\,\,\,g_{t_{i}}(x,z_{1},\ldots,z_{i-1},0,z_{i+1},\ldots,z_{n})=0$\foreignlanguage{british}{
for any $(z_{1},\ldots,z_{i-1},z_{i+1},\ldots,z_{n})\in(\mathbb{R}_{+})^{n-1}$;}\medskip{}

\noindent \textup{(iii)}$\,\,\,D_{\xi}f(x,0)=0$\foreignlanguage{british}{;}\medskip{}

\noindent \textup{(iv)}$\,\,\,g(x,t_{1},\ldots,t_{n})$\foreignlanguage{british}{
is convex with respect to $t_{i}\in\mathbb{R}_{+}$};\medskip{}

\noindent \textup{(v)}$\,\,\,g(x,t_{1},\ldots,t_{n})$\foreignlanguage{british}{
is non-decreasing with respect to $t_{i}\in\mathbb{R}_{+}$;}\medskip{}

\noindent \textup{(vi)}$\,\,\,\rho\mapsto f(x,\rho\,\xi)$\foreignlanguage{british}{
is non-decreasing in $\mathbb{R}_{+}$ for every $\xi\in\mathbb{R}^{m\times n}$.}
\end{singlespace}
\end{lem}

\begin{singlespace}
\noindent \begin{proof}[\bfseries{Proof}]We first prove (i), (ii)
and (iii). Let $\{\mathsf{e}_{\alpha}\}_{1\,\leq\,\alpha\,\leq m}$
be the standard basis of $\mathbb{R}^{m}$. Fix \foreignlanguage{english}{$\alpha\in\{1,...,m\}$}
and let $\tilde{g}:\Omega\times\mathbb{R}^{n}\to\mathbb{R}_{+}$ be
the function defined by 
\begin{equation}
\tilde{g}(x,t_{1},t_{2},\ldots,t_{n}):=f\left(x,\left(t_{1}\,\mathsf{e}_{\alpha}^{T},t_{2}\,\mathsf{e}_{\alpha}^{T},\ldots,t_{n}\,\mathsf{e}_{\alpha}^{T}\right)\right),\,\,\,\,\,\,\,\,\,\,\,\,x\in\Omega,\,\,(t_{1},...,t_{n})\in\mathbb{R}^{n}.\label{eq:gtilde}
\end{equation}
Then, for a.e. $x\in\Omega$, the function $\tilde{g}(x,\cdot)$ is
of class $C^{1}$ on $\mathbb{R}^{n}$, since it is a composition
of $C^{1}$ functions. Moreover, for every $i\in\{1,\ldots,n\}$ and
every $(z_{1},\ldots,z_{n})\in\mathbb{R}^{n}$, we get 
\begin{equation}
\tilde{g}_{t_{i}}(x,z_{1},\ldots,z_{i},\ldots,z_{n})\,=\,\frac{\partial f}{\partial\xi_{i}^{\alpha}}\left(x,\left(z_{1}\,\mathsf{e}_{\alpha}^{T},\ldots,z_{i}\,\mathsf{e}_{\alpha}^{T},\ldots,z_{n}\,\mathsf{e}_{\alpha}^{T}\right)\right).\label{eq:gtilde_derivate}
\end{equation}
Now observe that, for a.e. $x\in\Omega$, the partial map 
\[
(t_{1},...,t_{n})\in(\mathbb{R}_{+})^{n}\longmapsto g(x,t_{1},\ldots,t_{n})
\]
is the restriction of $\tilde{g}(x,\cdot)$ to $(\mathbb{R}_{+})^{n}$.
Indeed, by $(\mathbf{A2})$ and (\ref{eq:gtilde}),  for a.e. $x\in\Omega$
and every $(t_{1},...,t_{n})\in(\mathbb{R}_{+})^{n}$ we have
\[
g(x,t_{1},...,t_{n})\,=\,f\left(x,\left(t_{1}\,\mathsf{e}_{\alpha}^{T},\ldots,t_{n}\,\mathsf{e}_{\alpha}^{T}\right)\right)=\,\tilde{g}(x,t_{1},\ldots,t_{n})\,.
\]
Therefore, for a.e. $x\in\Omega$, the map $g(x,\cdot)$ is of class
$C^{1}$ on $(\mathbb{R}_{+})^{n}$.\\
Now also fix $i\in\{1,\ldots,n\}$, $(\xi_{1},\ldots,\xi_{i-1},\xi_{i+1},\ldots,\xi_{n})\in\mathbb{R}^{m\times(n-1)}$
and $x\in\Omega$, and consider the function $F:\mathbb{R}\rightarrow\mathbb{R}_{+}$
defined by
\[
F(t):=f(x,(\xi_{1},\ldots,\xi_{i-1},t\,\mathsf{e}_{\alpha}^{T},\xi_{i+1},\ldots,\xi_{n}))\,,\,\,\,\,\,\,\,\,\,\,t\in\mathbb{R},
\]
with the usual modications in the limit cases $i\in\{1,n\}$. Then
$F$ is of class $C^{1}$ on $\mathbb{R}$, because it is a composition
of $C^{1}$ functions. Furthermore, $F$ is an even function, since\begin{align*}
F(t)&=g(x,\vert\xi_{1}\vert,\ldots,\vert\xi_{i-1}\vert,\vert t\vert,\vert\xi_{i+1}\vert,\ldots,\vert\xi_{n}\vert)\\
&=g(x,\vert\xi_{1}\vert,\ldots,\vert\xi_{i-1}\vert,\vert-t\vert,\vert\xi_{i+1}\vert,\ldots,\vert\xi_{n}\vert)=F(-t)\,\,\,\,\,\,\,\,\,\,\mathrm{for\,\,every}\,\,t\in\mathbb{R}.
\end{align*}Hence,
\begin{equation}
F'(0)\,=\,\frac{\partial f}{\partial\xi_{i}^{\alpha}}(x,(\xi_{1},\ldots,\xi_{i-1},\underbrace{0^{T}}_{\in\,\mathbb{R}^{m}},\xi_{i+1},\ldots,\xi_{n}))=0\,.\label{eq:grad_f}
\end{equation}
Since $(\xi_{1},\ldots,\xi_{i-1},\xi_{i+1},\ldots,\xi_{n})\in\mathbb{R}^{m\times(n-1)}$
is arbitrary, from (\ref{eq:grad_f}) and (\ref{eq:gtilde_derivate})
with $z_{i}=0$ we deduce 
\[
g_{t_{i}}(x,z_{1},\ldots,z_{i-1},0,z_{i+1},\ldots,z_{n})\,=\,\frac{\partial f}{\partial\xi_{i}^{\alpha}}\left(x,\left(z_{1}\,\mathsf{e}_{\alpha}^{T},\ldots,z_{i-1}\,\mathsf{e}_{\alpha}^{T},0^{T},z_{i+1}\,\mathsf{e}_{\alpha}^{T},\ldots,z_{n}\,\mathsf{e}_{\alpha}^{T}\right)\right)=0
\]
for every $(z_{1},\ldots,z_{i-1},z_{i+1},\ldots,z_{n})\in(\mathbb{R}_{+})^{n-1}$,
where we have used the fact that $g(x,\cdot)$ is the restriction
of $\tilde{g}(x,\cdot)$ to $(\mathbb{R}_{+})^{n}$. The conclusions
(ii) and (iii) then follow from the arbitrariness of $x\in\Omega$,
$i\in\{1,\ldots,n\}$, $\alpha\in\{1,...,m\}$ and $(z_{1},\ldots,z_{i-1},z_{i+1},\ldots,z_{n})\in(\mathbb{R}_{+})^{n-1}$.\\
$\hspace*{1em}$We now prove (iv) and (v). Let $w\in\mathbb{R}^{m}$
with $\vert w\vert=1$ and fix $i\in\{1,\ldots,n\}$. Moreover, let
\[
(t_{1},\ldots,t_{i-1},a,b,t_{i+1},\ldots,t_{n})\,\in\,(\mathbb{R}_{+})^{n+1},
\]
\[
\zeta=(t_{1}\,w^{T},\ldots,t_{i-1}\,w^{T},a\,w^{T},t_{i+1}\,w^{T},\ldots,t_{n}\,w^{T})\,,
\]
\[
\eta=(t_{1}\,w^{T},\ldots,t_{i-1}\,w^{T},b\,w^{T},t_{i+1}\,w^{T},\ldots,t_{n}\,w^{T})\,.
\]
Then, for a.e. $x\in\Omega$ and for all $\theta\in[0,1]$ we have\begin{align*}
&g(x,t_{1},\ldots,t_{i-1},\theta a+(1-\theta)b,t_{i+1},\ldots,t_{n})\\
&\,\,\,\,\,\,\,=\,f(x,\theta\,\zeta+(1-\theta)\,\eta)\\
&\,\,\,\,\,\,\,\le\,\theta\,f(x,\zeta)+(1-\theta)\,f(x,\eta)\\
&\,\,\,\,\,\,\,=\,\theta\,g(x,t_{1},\ldots,t_{i-1},a,t_{i+1},\ldots,t_{n})\,+\,(1-\theta)\,g(x,t_{1},\ldots,t_{i-1},b,t_{i+1},\ldots,t_{n})\,.
\end{align*} This implies that $g(x,t_{1},\ldots,t_{n})$ is convex with respect
to each variable $t_{i}\in\mathbb{R}_{+}$. Therefore, the partial
maps 
\[
z\in\mathbb{R}_{+}\longmapsto g_{t_{i}}(x,t_{1},\ldots,t_{i-1},z,t_{i+1},\ldots,t_{n})\,,\,\,\,\,\,\,\,\,\,\,i\in\{1,\ldots,n\},
\]
are non-decreasing and this monotonicity property, together with (i)
and (ii), entails that $g(x,t_{1},\ldots,t_{n})$ is non-decreasing
with respect to each variable $t_{i}\in\mathbb{R}_{+}$.\\
$\hspace*{1em}$Finally, let us prove (vi). Fix $\rho_{1},\rho_{2}\in\mathbb{R}_{+}$
with $\rho_{1}<\rho_{2}$. Then, for a.e. $x\in\Omega$ and every
$\xi\in\mathbb{R}^{m\times n}$ we have\begin{align*}
f(x,\rho_{1}\,\xi)\,&=\,g(x,\rho_{1}\vert\xi_{1}\vert,\rho_{1}\vert\xi_{2}\vert,\ldots,\rho_{1}\vert\xi_{n-1}\vert,\rho_{1}\vert\xi_{n}\vert)\\
&\le\,g(x,\rho_{2}\vert\xi_{1}\vert,\rho_{1}\vert\xi_{2}\vert,\ldots,\rho_{1}\vert\xi_{n-1}\vert,\rho_{1}\vert\xi_{n}\vert)\\
&\leq\,\cdots\,\leq\,g(x,\rho_{2}\vert\xi_{1}\vert,\rho_{2}\vert\xi_{2}\vert,\ldots,\rho_{2}\vert\xi_{n-1}\vert,\rho_{1}\vert\xi_{n}\vert)\\
&\leq\,g(x,\rho_{2}\vert\xi_{1}\vert,\rho_{2}\vert\xi_{2}\vert,\ldots,\rho_{2}\vert\xi_{n-1}\vert,\rho_{2}\vert\xi_{n}\vert)\,=\,f(x,\rho_{2}\,\xi)\,,
\end{align*}where, in the last three lines, we have repeatedly used the property
\foreignlanguage{english}{(v)}. This concludes the proof.\end{proof}

\noindent $\hspace*{1em}$\foreignlanguage{english}{If $f$ is as
in Section \ref{assumptions}, then $W^{1,\mathcal{F}}(\Omega;\mathbb{R}^{m})$
is a vector space; this is a consequence of the following lemma.}
\end{singlespace}
\selectlanguage{english}%
\begin{lem}
\begin{singlespace}
\noindent \label{lemma1:daipotesi} Assume $(\mathbf{A1})$, $(\mathbf{A2})$
and $(\mathbf{A3})$. Then, for a.e. $x\in\Omega$\foreignlanguage{british}{,}\medskip{}

\noindent \textup{(i)}$\,\,\,f(x,\gamma\xi)\le\max\,\{1,\gamma^{\tau}\}f(x,\xi)$
for every $\gamma>0$ and every $\xi\in\mathbb{R}^{m\times n}$;\medskip{}

\noindent \textup{(ii)}$\,\,\,f(x,\xi+\eta)\le2^{\tau-1}\left[f(x,\xi)+f(x,\eta)\right]$
for every $\xi,\eta\in\mathbb{R}^{m\times n}$,\medskip{}

\noindent where $\tau\geq1$ is the constant in $(\mathbf{A3})$.
\end{singlespace}
\end{lem}

\selectlanguage{british}%
\noindent \begin{proof}[\bfseries{Proof}] Let us prove (i). If $\xi=0$,
then $f(x,\gamma\xi)=f(x,\xi)$ for every \foreignlanguage{english}{$\gamma>0$}
and the conclusion immediately follows.\\
Assume $\vert\xi\vert>0$. We consider the cases $\gamma>1$ and $0<\gamma\leq1$
separately.\\
$\hspace*{1em}$Let $\gamma>1$. Then $(\mathbf{A3})$ implies $f(x,\gamma\xi)\leq\gamma^{\tau}f(x,\xi)$.\\
If instead $0<\gamma\leq1$, by Lemma \ref{lem:g function} (vi) we
get $f(x,\gamma\xi)\leq f(x,\xi)$ and the conclusion follows.\\
$\hspace*{1em}$Let us now prove (ii). If $\xi,\eta\in\mathbb{R}^{m\times n}$,
by \foreignlanguage{english}{$(\mathbf{A1})$} and \foreignlanguage{english}{$(\mathbf{A3})$}
we obtain 
\[
f(x,\xi+\eta)\leq\frac{1}{2}\left[f(x,2\xi)+f(x,2\eta)\right]\leq2^{\tau-1}\left[f(x,\xi)+f(x,\eta)\right].
\]
This completes the proof.\end{proof}

\begin{singlespace}
\noindent $\hspace*{1em}$\foreignlanguage{english}{We now recall
the following elementary result, whose proof can be found, for example,
in \cite[Lemma 3.1]{cupmarmas1}.}
\end{singlespace}
\selectlanguage{english}%
\begin{lem}
\begin{singlespace}
\noindent \label{delta2nuovo} Consider $h:\mathbb{R}_{+}\to\mathbb{R}_{+}$
of class $C^{1}$. Assume that there exists $\tau\geq1$ such that
\[
h(\gamma t)\le\gamma^{\tau}h(t)\,\,\,\,\,\,\,\mathit{for\,\,all\,\,}\gamma>1\,\,\mathit{and}\,\,t\geq0.
\]
Then 
\[
h^{\prime}(t)\,t\le\tau\,h(t)\,\,\,\,\,\,\,\mathit{for\,\,all\,\,}t\geq0.
\]
\end{singlespace}
\end{lem}

\selectlanguage{british}%

\subsection{Anisotropic Sobolev spaces}

$\hspace*{1em}$Let $q_{i}\geq1$ for all $i\in\{1,\ldots,n\}$. For
any open subset $\Omega$ of $\mathbb{R}^{n}$, we consider the anisotropic
Sobolev space 
\[
W^{1,(q_{1},\ldots,q_{n})}(\Omega;\mathbb{R}^{m}):=\left\{ v\in W^{1,1}(\Omega;\mathbb{R}^{m}):v_{x_{i}}\in L^{q_{i}}(\Omega;\mathbb{R}^{m}),\mathrm{\,\,for\,\,all\,\,}i=1,\ldots,n\right\} 
\]
endowed with the natural norm
\[
\Vert v\Vert_{W^{1,(q_{1},\ldots,q_{n})}(\Omega;\mathbb{R}^{m})}:=\Vert v\Vert_{L^{1}(\Omega;\mathbb{R}^{m})}+\sum_{i=1}^{n}\Vert v_{x_{i}}\Vert_{L^{q_{i}}(\Omega;\mathbb{R}^{m})}\,.
\]
Sometimes, when no confusion may arise, we will omit the target space
$\mathbb{R}^{m}$. Let us denote $\mathbf{q}=(q_{1},...,q_{n})$ and
\[
W_{0}^{1,(q_{1},\ldots,q_{n})}(\Omega;\mathbb{R}^{m})=W_{0}^{1,1}(\Omega;\mathbb{R}^{m})\cap W^{1,(q_{1},\ldots,q_{n})}(\Omega;\mathbb{R}^{m})\,.
\]
These spaces are studied in \cite{Troisi} (see also \cite{AceFus}).
We now report an embedding theorem for this class of spaces, whose
proof can be obtained by a straightforward adaptation of that of \cite[Theorem 1.2]{Troisi}. 
\begin{thm}
\label{thm:embedding1}Let $\Omega\subset\mathbb{R}^{n}$ be a bounded
open set and let $v\in W_{0}^{1,(q_{1},\ldots,q_{n})}(\Omega;\mathbb{R}^{m})$,
$q_{i}\geq1$ for all $i\in\{1,\ldots,n\}$. Then,\medskip{}

\noindent $\mathrm{(}1\mathrm{)}\,$ if $\overline{\mathbf{q}}<n$,
we have 
\[
\Vert v\Vert_{L^{\overline{\mathbf{q}}^{*}}(\Omega;\mathbb{R}^{m})}\leq\,c\,\prod_{i=1}^{n}\Vert v_{x_{i}}\Vert_{L^{q_{i}}(\Omega;\mathbb{R}^{m})}^{\frac{1}{n}}
\]
for a constant $c=c(m,n,\mathbf{q})>0$, where $\overline{\mathbf{q}}^{*}$
is defined by $(\ref{eq:SobConj})$$-$$(\ref{eq:harmonic})$ with
$\mathbf{p}=\mathbf{q}$ and $\ell=\overline{\mathbf{q}}$;\medskip{}

\begin{singlespace}
\noindent $\mathrm{(}2\mathrm{)}\,$ if $\overline{\mathbf{q}}\geq n$,
for every $1\leq\chi<\infty$ we have 
\[
\Vert v\Vert_{L^{\chi}(\Omega;\mathbb{R}^{m})}\leq\,c\,\vert\Omega\vert^{\frac{1}{\chi}\,+\,\frac{1}{n}\,-\,\frac{1}{\overline{\mathbf{q}}}}\,\prod_{i=1}^{n}\Vert v_{x_{i}}\Vert_{L^{q_{i}}(\Omega;\mathbb{R}^{m})}^{\frac{1}{n}}
\]
for a constant $c=c(m,n,\mathbf{q},\chi)>0$.
\end{singlespace}
\end{thm}

\noindent $\hspace*{1em}$The following embedding result is proved
in \cite{AceFus}.
\begin{thm}
\noindent \label{thm:embedding2}Let $Q\subset\mathbb{R}^{n}$ be
a cube with edges parallel to the coordinate axes and consider a function
$v\in W^{1,(q_{1},\ldots,q_{n})}(Q;\mathbb{R}^{m})$, $q_{i}\geq1$
for all $i\in\{1,\ldots,n\}$. Let $\max\,\{q_{i}\}<\overline{\mathbf{q}}^{*}$.
Then $v\in L^{\overline{\mathbf{q}}^{*}}(Q;\mathbb{R}^{m})$. Moreover,
there exists a positive constant $c$ depending on $n$ and $\mathbf{q}$
if $\overline{\mathbf{q}}<n$, and additionally on $Q$ if $\overline{\mathbf{q}}\geq n$,
such that 
\[
\Vert v\Vert_{L^{\overline{\mathbf{q}}^{*}}(Q)}\leq\,c\left\{ \Vert v\Vert_{L^{1}(Q)}+\sum_{i=1}^{n}\Vert v_{x_{i}}\Vert_{L^{q_{i}}(Q)}\right\} .
\]
\end{thm}

\noindent $\hspace*{1em}$The next proposition is a consequence of
the above theorem.
\selectlanguage{english}%
\begin{prop}
\begin{singlespace}
\noindent \label{prop:summability}\foreignlanguage{british}{For every
$i\in\{1,\ldots,n\}$, let $p_{i}\in(1,\infty)$, $r_{i}\in[1,\infty]$,
with $r_{i}\geq\frac{1}{p_{i}-1}$ if $1<p_{i}<2$, and let $\sigma_{i}$
be defined according to $(\ref{eq:sigma_i})$.} Moreover, assume that
\foreignlanguage{british}{$\max\,\{\sigma_{i}\}<\overline{\sigma}^{*}$}
and that the left inequality in $(\ref{growth})$ holds under the
summability conditions
\[
\lambda_{i}^{-1}\in L_{loc}^{r_{i}}(\Omega),\,\,\,\,\,{\color{red}{\normalcolor i\in\{1,...,n\}}}.
\]
If $u\in W^{1,\mathcal{F}}(\Omega;\mathbb{R}^{m})$, then we get $\vert u\vert\in L_{loc}^{\overline{\sigma}^{*}}(\Omega)$.
\end{singlespace}
\end{prop}

\selectlanguage{british}%
\noindent \begin{proof}[\bfseries{Proof}] Let $Q\Subset\Omega$ be
a cube with edges parallel to the coordinate axes and fix $j\in\{1,\ldots,n\}$.\\
$\hspace*{1em}$Let us first assume that $r_{j}<\infty$. Recalling
that $\sigma_{j}=\frac{p_{j}r_{j}}{r_{j}+1}$, by Hölder's inequality
we obtain
\begin{equation}
\Vert u_{x_{j}}\Vert_{L^{\sigma_{j}}(Q)}^{p_{j}}=\left(\int_{Q}\lambda_{j}^{\frac{r_{j}}{r_{j}+1}}\,\vert u_{x_{j}}\vert^{\frac{p_{j}r_{j}}{r_{j}+1}}\,\lambda_{j}^{-\,\frac{r_{j}}{r_{j}+1}}dx\right)^{\frac{r_{j}+1}{r_{j}}}\leq\,\Vert\lambda_{j}^{-1}\Vert_{L^{r_{j}}(Q)}\int_{Q}\lambda_{j}\,\vert u_{x_{j}}\vert^{p_{j}}\,dx\,.\label{eq:Elv1}
\end{equation}
Combining the previous estimate with the left inequality in (\ref{growth}),
we have
\[
\Vert u_{x_{j}}\Vert_{L^{\sigma_{j}}(Q)}^{p_{j}}\leq\,\Vert\lambda_{j}^{-1}\Vert_{L^{r_{j}}(Q)}\,\sum_{i=1}^{n}\int_{\Omega}\lambda_{i}\,\vert u_{x_{i}}\vert^{p_{i}}\,dx\,\leq\,\Vert\lambda_{j}^{-1}\Vert_{L^{r_{j}}(Q)}\int_{\Omega}f(x,Du)\,dx\,,
\]
and the last term is finite, because $\lambda_{j}^{-1}\in L_{loc}^{r_{j}}(\Omega)$
and \foreignlanguage{english}{$u\in W^{1,\mathcal{F}}(\Omega;\mathbb{R}^{m})$}.\\
$\hspace*{1em}$If $r_{j}=\infty$, then $\sigma_{j}=p_{j}$ and we
find
\begin{equation}
\Vert u_{x_{j}}\Vert_{L^{p_{j}}(Q)}^{p_{j}}=\int_{Q}\lambda_{j}\,\vert u_{x_{j}}\vert^{p_{j}}\,\lambda_{j}^{-1}\,dx\,\leq\,\Vert\lambda_{j}^{-1}\Vert_{L^{\infty}(Q)}\int_{Q}\lambda_{j}\,\vert u_{x_{j}}\vert^{p_{j}}\,dx\,,\label{eq:Elv2}
\end{equation}
which combined with the left inequality in (\ref{growth}) gives
\[
\Vert u_{x_{j}}\Vert_{L^{p_{j}}(Q)}^{p_{j}}\leq\,\Vert\lambda_{j}^{-1}\Vert_{L^{\infty}(Q)}\int_{\Omega}f(x,Du)\,dx\,<+\infty\,.
\]
We have thus proved that $u\in W^{1,(\sigma_{1},\ldots,\sigma_{n})}(Q;\mathbb{R}^{m})$.
Since $\sigma_{i}\geq1$ for all $i\in\{1,\ldots,n\}$ and $\max\,\{\sigma_{i}\}<\overline{\sigma}^{*}$,
by Theorem \ref{thm:embedding2} we get $u\in L^{\overline{\mathbf{\sigma}}^{*}}(Q;\mathbb{R}^{m})$.
The conclusion then follows from the arbitrariness of $Q$.\end{proof}

\noindent $\hspace*{1em}$We conclude this section with the following
proposition about a Poincaré-Sobolev type inequality.
\begin{prop}
\noindent \label{prop:Poincar=0000E9}For every $i\in\{1,\ldots,n\}$,
let $p_{i}\in(1,\infty)$, $r_{i}\in[1,\infty]$, with $r_{i}\geq\frac{1}{p_{i}-1}$
if $1<p_{i}<2$, and let $\sigma_{i}$ be defined according to $(\ref{eq:sigma_i})$.
Consider a bounded open set $\Omega\subset\mathbb{R}^{n}$ and let
$v\in W_{0}^{1,(\sigma_{1},\ldots,\sigma_{n})}(\Omega;\mathbb{R}^{m})$.
Moreover, for every $i\in\{1,\ldots,n\}$, let \foreignlanguage{english}{$\lambda_{i}:\Omega\rightarrow[0,\infty)$
be a measurable function such that $\lambda_{i}^{-1}\in L^{r_{i}}(\Omega)$.
Then, there exists }a positive constant $c$, depending on $m$, $n$
and $\sigma=(\sigma_{1},\ldots,\sigma_{n})$, such that\foreignlanguage{english}{
\[
\left(\int_{\Omega}\vert v\vert^{\overline{\sigma}^{*}}dx\right)^{\frac{1}{\overline{\sigma}^{*}}}\leq\,c\,\mathfrak{M}\,\prod_{i=1}^{n}\Vert\lambda_{i}^{-1}\Vert_{L^{r_{i}}(\Omega)}^{\frac{1}{np_{i}}}\left(\int_{\Omega}\lambda_{i}\,\vert v_{x_{i}}\vert^{p_{i}}\,dx\right)^{\frac{1}{np_{i}}},
\]
where 
\[
\mathfrak{M}:=\begin{cases}
\begin{array}{cc}
1 & \mathit{if}\,\,\overline{\sigma}<n,\\
\vert\Omega\vert^{\frac{1}{\overline{\sigma}^{*}}\,+\,\frac{1}{n}\,-\,\frac{1}{\overline{\mathbf{\sigma}}}} & \mathit{if}\,\,\overline{\sigma}\geq n.
\end{array}\end{cases}
\]
}
\end{prop}

\noindent \begin{proof}[\bfseries{Proof}]Notice that the assumptions
imply $\sigma_{i}\geq1$ for all $i\in\{1,\ldots,n\}$. Therefore,
by Theorem \ref{thm:embedding1} we get
\[
\Vert v\Vert_{L^{\overline{\sigma}^{*}}(\Omega)}\leq\,c\,\mathfrak{M}\,\prod_{i=1}^{n}\Vert v_{x_{i}}\Vert_{L^{\sigma_{i}}(\Omega)}^{\frac{1}{n}}
\]
for a positive constant $c$ depending on $m$, $n$ and $\sigma$.
Arguing as in (\ref{eq:Elv1})\foreignlanguage{english}{$-$}(\ref{eq:Elv2}),
we obtain
\[
\Vert v\Vert_{L^{\overline{\sigma}^{*}}(\Omega)}\leq\,c\,\mathfrak{M}\,\prod_{i=1}^{n}\Vert v_{x_{i}}\Vert_{L^{\sigma_{i}}(\Omega)}^{\frac{1}{n}}\,\leq\,c\,\mathfrak{M}\,\prod_{i=1}^{n}\Vert\lambda_{i}^{-1}\Vert_{L^{r_{i}}(\Omega)}^{\frac{1}{np_{i}}}\left(\int_{\Omega}\lambda_{i}\,\vert v_{x_{i}}\vert^{p_{i}}\,dx\right)^{\frac{1}{np_{i}}}.
\]
This concludes the proof.\end{proof}
\selectlanguage{english}%
\begin{singlespace}

\section{The Euler's equation\label{s:eulero}}
\end{singlespace}

\selectlanguage{british}%
\begin{singlespace}
\noindent $\hspace*{1em}$In this section, we prove that an Euler's
equation holds true. This equation will be our starting point in the
proof of Theorem \ref{thm:theo1}.\medskip{}

\end{singlespace}
\selectlanguage{english}%
\begin{prop}
\begin{singlespace}
\noindent \label{Euler'} Assume $(\mathbf{A1})$$-$$(\mathbf{A3})$
and let $u$ be a local minimizer of $(\ref{functional})$. Then 
\begin{eqnarray}
 &  & \int_{\Omega}\sum_{i=1}^{n}\sum_{\alpha=1}^{m}\frac{\partial f}{\partial{\xi_{i}^{\alpha}}}\left(x,Du\right)\,\varphi_{x_{i}}^{\alpha}\,dx=0\label{p:eulero'}
\end{eqnarray}
for all $\varphi\in W^{1,\mathcal{F}}(\Omega;\mathbb{R}^{m})$ with
${\color{red}{\normalcolor \mathrm{supp}\,\varphi\Subset\Omega}}$. 
\end{singlespace}
\end{prop}

\selectlanguage{british}%
\begin{singlespace}
\noindent \begin{proof}[\bfseries{Proof}]\foreignlanguage{english}{Let
$\varphi\in W^{1,\mathcal{F}}(\Omega;\mathbb{R}^{m})$ with ${\color{red}{\normalcolor \mathrm{supp}\,\varphi\Subset\Omega}}$.
By (\ref{eq:struttura}) also $-\varphi$ is in $W^{1,\mathcal{F}}(\Omega;\mathbb{R}^{m})$.
By Lemma \ref{lemma1:daipotesi} we get $u+t\varphi\in W^{1,\mathcal{F}}(\Omega;\mathbb{R}^{m})$
for every $t\in\mathbb{R}$. By the local minimality of $u$, 
\[
\mathcal{F}(u)\,\leq\,\mathcal{F}(u+t\varphi)\qquad\forall\,t\in\mathbb{R}.
\]
To prove \eqref{p:eulero'} it suffices to prove that 
\[
\left.\frac{d}{dt}\,\mathcal{F}(u+t\varphi)\right\vert _{t=0}=\int_{\Omega}\left.\frac{d}{dt}f(x,Du(x)+t\,D\varphi(x))\right\vert _{t=0}\,dx\,.
\]
Therefore, we need to prove that 
\[
\Big\vert\sum_{i=1}^{n}\sum_{\alpha=1}^{m}\frac{\partial f}{\partial{\xi_{i}^{\alpha}}}(x,Du+t\,D\varphi)\varphi_{x_{i}}^{\alpha}\Big\vert\leq H(x)\qquad\forall\,t\in(-1,1)
\]
with $H\in L^{1}(\Omega)$. By the convexity of $f(x,\cdot)$, we
obtain 
\[
f(x,\xi_{0})-f(x,2\xi_{0}-\xi)\leq\sum_{i=1}^{n}\sum_{\alpha=1}^{m}\frac{\partial f}{\partial{\xi_{i}^{\alpha}}}(x,\xi_{0})(\xi_{i}^{\alpha}-(\xi_{0})_{i}^{\alpha})\leq f(x,\xi)-f(x,\xi_{0})\,.
\]
If $\xi_{0}=Du(x)+t\,D\varphi(x)$ and $\xi=Du(x)+(1+t)\,D\varphi(x)$,
we have 
\[
2\xi_{0}-\xi=Du(x)+(t-1)\,D\varphi(x)
\]
and\begin{align*}
f(x,Du+t\,D\varphi)-f(x,Du+(t-1)\,D\varphi)& \leq \sum_{i=1}^{n}\sum_{\alpha=1}^{m}\frac{\partial f}{\partial{\xi_{i}^{\alpha}}}(x,Du+t\,D\varphi)\varphi_{x_{i}}^{\alpha}\\
& \leq f(x,Du+(1+t)\,D\varphi)-f(x,Du+t\,D\varphi)\,.
\end{align*}Therefore, since $f$ is non-negative,\vspace{1mm}}

\selectlanguage{english}%
\noindent 
\[
\Big|\sum_{i=1}^{n}\sum_{\alpha=1}^{m}\frac{\partial f}{\partial{\xi_{i}^{\alpha}}}(x,Du+t\,D\varphi)\varphi_{x_{i}}^{\alpha}\Big|\leq f(x,Du+(1+t)\,D\varphi)+f(x,Du+(t-1)\,D\varphi)\,.
\]
\foreignlanguage{british}{Moreover, since $t\in(-1,1)$, we can use
Lemma} \ref{lemma1:daipotesi} \foreignlanguage{british}{to estimate
the first term on the right-hand side of the previous inequality,
thus obtaining\begin{align*}
f(x,Du+(1+t)\,D\varphi)&\leq2^{\tau-1}\left[f(x,Du)+f(x,(1+t)\,D\varphi)\right]\\
&\leq2^{\tau-1}\left[f(x,Du)+2^{\tau}f(x,D\varphi)\right].
\end{align*}To estimate $f(x,Du\,+\,(t-1)\,D\varphi)$, we now consider the cases
$t\in[0,1)$ and} \foreignlanguage{british}{$t\in(-1,0)$ separately.
Let us first assume $t\in[0,1)$. Then},\foreignlanguage{british}{
}using again the convexity of $f(x,\cdot)$ and Lemma \ref{lemma1:daipotesi},
we get\vspace{-2mm}\begin{align*}
f(x,Du+(t-1)\,D\varphi)&\leq t\,f(x,Du)+(1-t)f(x,Du-D\varphi)\\
&\leq f(x,Du)+f(x,Du-D\varphi)\\
&\leq (1+2^{\tau-1})f(x,Du)+2^{\tau-1}f(x,-D\varphi)\,.
\end{align*}Now assume \foreignlanguage{british}{$t\in(-1,0)$. Since $1<1-t<2$,
by Lemma} \ref{lemma1:daipotesi} we have\begin{align*}
f(x,Du+(t-1)\,D\varphi) & \leq 2^{\tau-1}\left[f(x,Du)+f(x,(t-1)\,D\varphi)\right]\\
&\leq2^{\tau-1}\left[f(x,Du)+2^{\tau}f(x,-D\varphi)\right]\\
&= 2^{\tau-1}f(x,Du)+2^{2\tau-1}f(x,-D\varphi)\,.
\end{align*}We have so proved that
\[
\Big\vert\sum_{i=1}^{n}\sum_{\alpha=1}^{m}\frac{\partial f}{\partial{\xi_{i}^{\alpha}}}(x,Du+t\,D\varphi)\varphi_{x_{i}}^{\alpha}\Big\vert\leq(1+2^{\tau})f(x,Du)+2^{2\tau-1}[f(x,D\varphi)+f(x,-D\varphi)]=:H(x)\,.
\]
Since $u,\varphi,-\varphi\in W^{1,\mathcal{F}}(\Omega;\mathbb{R}^{m})$,
we conclude\foreignlanguage{british}{.\end{proof}}
\end{singlespace}
\selectlanguage{british}%
\begin{singlespace}

\section{Proof of the main result\label{sec:Proof}}
\end{singlespace}

\begin{singlespace}
\noindent $\hspace*{1em}$We first state a lemma useful for the proof
of Theorem \ref{thm:theo1}. In the statement, the functions $\lambda_{i},\mu$
and the exponents \foreignlanguage{english}{$\{p_{i}\},q,\{r_{i}\},s$
}are the same as in the statement of Theorem \ref{thm:theo1}. Moreover,
$x_{0}\in\Omega$ and $R_{0}\in(0,1]$ are such that $B_{R_{0}}:=B_{R_{0}}(x_{0})\Subset\Omega$.
Fixed $0<\rho<R\le R_{0}$, the function $\eta\in C_{c}^{\infty}(B_{R}(x_{0}))$
denotes a cut-off function satisfying \foreignlanguage{english}{
\begin{equation}
0\leq\eta\leq1,\qquad\eta\equiv1\,\,\,\,\text{on \ensuremath{B_{\rho}(x_{0})},}\qquad\text{\ensuremath{|D\eta|\leq\frac{2}{R-\rho}}}\,.\label{eq:eta}
\end{equation}
}
\end{singlespace}
\selectlanguage{english}%
\begin{lem}
\begin{singlespace}
\noindent \label{lem:lemma prelim}Let us assume that $(\mathbf{A1})$$-$$(\mathbf{A4})$
hold under the summability conditions
\[
\lambda_{i}^{-1}\in L_{loc}^{r_{i}}(\Omega),\,\,\,{\color{red}{\normalcolor i\in\{1,...,n\}}},\;\;\;\;\;\;{\color{red}{\normalcolor \mu\in L_{loc}^{s}(\Omega)}},
\]
for some ${\color{red}{\normalcolor r_{i}\in[1,\infty]}}$ and $s\in(1,\infty]$
such that 
\[
\max\,\{\sigma_{i}\}<\overline{\sigma}^{*}\,\,\,\,\,\,\,\,\,\,\mathit{and}\,\,\,\,\,\,\,\,\,\,\frac{1}{\overline{\mathbf{pr}}}\,+\,\frac{1}{qs}\,+\,\frac{1}{\overline{\mathbf{p}}}\,-\,\frac{1}{q}\,<\,\frac{1}{n}\,.
\]
\foreignlanguage{british}{For every }$i\in\{1,\ldots,n\}$,\foreignlanguage{british}{
}if\foreignlanguage{british}{ $1<p_{i}<2$} we also require\foreignlanguage{british}{
$r_{i}\geq\frac{1}{p_{i}-1}$}. Moreover, let $u$ be a local minimizer
of $(\ref{functional})$. Then, for every $\gamma\geq0$ and every
$i\in\{1,\ldots,n\}$ we have\foreignlanguage{british}{
\begin{equation}
\int_{B_{R}}\lambda_{i}(x)\,\vert u\vert^{p_{i}\gamma}\,\vert u_{x_{i}}\vert^{p_{i}}\,\eta^{q}\,dx\,\leq\,\frac{c}{(R-\rho)^{q}}\int_{B_{R}}\mu(x)\left(\max\,\{\vert u\vert,1\}\right)^{q+p_{i}\gamma}dx\label{eq:jolly}
\end{equation}
for a positive constant $c$ depending on $n$, $q$ and $\tau$,
but independent of $\gamma$, $u$, $R$ and $\rho$.}
\end{singlespace}
\end{lem}

\selectlanguage{british}%
\begin{singlespace}
\noindent \begin{proof}[\bfseries{Proof}] We begin by defining a
class of suitable test functions for the Euler's equation (\ref{p:eulero'}).
Let us approximate the identity \foreignlanguage{english}{function
$\mathrm{id}:\mathbb{R}_{+}\rightarrow\mathbb{R}_{+}$ with an increasing
sequence of non-decreasing $C^{1}$ functions $h_{k}:\mathbb{R}_{+}\rightarrow\mathbb{R}_{+}$
having the following properties:
\[
h_{k}(t)=0\quad\forall\,t\in\left[0,\frac{1}{k}\right],\quad h_{k}(t)=k\quad\forall\,t\in[k+1,\infty),\quad0\leq h_{k}^{\prime}(t)\leq2\,\,\,\,\text{in}\ \mathbb{R}_{+}\,.
\]
}

\selectlanguage{english}%
\noindent Fixed $k,i\in\mathbb{N}$, $i\leq n$, and $\gamma>0$,
let $\Phi_{k}^{(i,\gamma)}:\mathbb{R}_{+}\rightarrow\mathbb{R}_{+}$
be the increasing function defined as follows
\[
\Phi_{k}^{(i,\gamma)}(t):=h_{k}(t^{p_{i}\gamma})\,.
\]
Define $\varphi_{k}^{(i,\gamma)}:B_{R}(x_{0})\rightarrow\mathbb{R}^{m}$
by
\[
\varphi_{k}^{(i,\gamma)}(x):=\Phi_{k}^{(i,\gamma)}(|u(x)|)\,u(x)\,[\eta(x)]^{q}.
\]
From now on, we omit the dependence of $\Phi_{k}$ and $\varphi_{k}$
on $i$ and $\gamma$, i.e. $\Phi_{k}=\Phi_{k}^{(i,\gamma)}$ and
$\varphi_{k}=\varphi_{k}^{(i,\gamma)}$. We have that $\Phi_{k}$
is in $C^{1}(\mathbb{R}_{+})$, bounded and with bounded derivative.
Precisely, define $a_{k}$ and $b_{k}$ positive real numbers, such
that $a_{k}^{p_{i}\gamma}=\frac{1}{k}$ and $b_{k}^{p_{i}\gamma}=k+1$.
In particular,
\[
\Phi'_{k}(s)=\begin{cases}
\begin{array}{cc}
0 & \mathrm{if}\,\,\,s\in\mathbb{R}_{+}\backslash\,[a_{k},b_{k}]\,,\\
p_{i}\gamma\,h'_{k}(s^{p_{i}\gamma})s^{p_{i}\gamma-1} & \mathrm{if}\,\,\,s\in[a_{k},b_{k}]\,,\,\,\,\,\,\,\,\,\,\,\,\,
\end{array}\end{cases}
\]
and 
\[
\Vert\Phi'_{k}\Vert_{L^{\infty}(\mathbb{R}_{+})}\le2p_{i}\gamma\,\max\,\{a_{k}^{p_{i}\gamma-1},b_{k}^{p_{i}\gamma-1}\}<+\infty\,.
\]
As a consequence, taking into account that $u\in W^{1,1}(\Omega;\mathbb{R}^{m})$,
we have that $\Phi_{k}(|u|)\,u$ is in $W^{1,1}(\Omega;\mathbb{R}^{m})$.
This implies that $\varphi_{k}\in W^{1,1}(\Omega;\mathbb{R}^{m})$
too. Moreover, ${\color{red}{\normalcolor \mathrm{supp}\,\varphi_{k}\Subset B_{R}(x_{0})}}.$
Therefore, by Lemma \ref{lemma1:daipotesi} (ii), \foreignlanguage{british}{$\varphi_{k}$}
\foreignlanguage{british}{is a test function for the Euler's equation
(\ref{p:eulero'})} if we prove that\begin{align*}
&I_{1}:=\int_{B_{R}\cap\{a_{k}\,<\,\vert u\vert\,<\,b_{k}\}}f\left(x,\Phi'_{k}(\vert u\vert)\,\frac{u(x)}{\vert u(x)\vert}\,\langle u,u_{x_{1}}\rangle\,\eta^{q},\ldots,\Phi'_{k}(\vert u\vert)\,\frac{u(x)}{\vert u(x)\vert}\,\langle u,u_{x_{n}}\rangle\,\eta^{q}\right)dx<+\infty\,,\\
&I_{2}:=\int_{B_{R}\cap\{\vert u\vert\,>\,a_{k}\}}f(x,\eta^{q}\,\Phi_{k}(\vert u\vert)\,Du)\,dx<+\infty\,,\\
&I_{3}:=\int_{B_{R}\cap\{\vert u\vert\,>\,a_{k}\}}f(x,\Phi_{k}(\vert u\vert)\,u\,q\eta^{q-1}\,\eta_{x_{1}},\ldots,\Phi_{k}(\vert u\vert)\,u\,q\eta^{q-1}\,\eta_{x_{n}})\,dx<+\infty\,.
\end{align*}Of course, by the definitions of $\Phi_{k}$, $a_{k}$, $b_{k}$ and
by (\ref{growth}), we have\begin{align*}
&\int_{B_{R}\cap\{\vert u\vert\,\notin\,(a_{k},b_{k})\}}f\left(x,\Phi'_{k}(\vert u\vert)\,\frac{u(x)}{\vert u(x)\vert}\,\langle u,u_{x_{1}}\rangle\,\eta^{q},\ldots,\Phi'_{k}(\vert u\vert)\,\frac{u(x)}{\vert u(x)\vert}\,\langle u,u_{x_{n}}\rangle\,\eta^{q}\right)dx\\
&\,\,\,\,\,\,\,=\int_{B_{R}}f(x,0)\,dx\,\leq\int_{B_{R}}\mu(x)\,dx\,<+\infty\,,\\
&\int_{B_{R}\cap\{\vert u\vert\,\leq\,a_{k}\}}f(x,\eta^{q}\,\Phi_{k}(\vert u\vert)\,Du)\,dx\,=\int_{B_{R}}f(x,0)\,dx\,<+\infty\,,\\
&\int_{B_{R}\cap\{\vert u\vert\,\leq\,a_{k}\}}f(x,\Phi_{k}(\vert u\vert)\,u\,q\eta^{q-1}\,\eta_{x_{1}},\ldots,\Phi_{k}(\vert u\vert)\,u\,q\eta^{q-1}\,\eta_{x_{n}})\,dx\,=\int_{B_{R}}f(x,0)\,dx\,<+\infty\,.
\end{align*}We now estimate $I_{1}$, $I_{2}$ and $I_{3}$ separately. Let us
first consider $I_{1}$. Using Lemma \ref{lemma1:daipotesi} (i),
$(\mathbf{A2})$, the Cauchy-Schwarz inequality, Lemma \ref{lem:g function}
(v) and $(\mathbf{A3})$, we find\begin{align*}
I_{1}&\le\max\,\{1,\Vert\Phi'_{k}\Vert_{L^{\infty}(\mathbb{R}_{+})}^{\tau}\}\int_{B_{R}\cap\{a_{k}\,<\,\vert u\vert\,<\,b_{k}\}}f\left(x,\frac{u}{\vert u\vert}\langle u,u_{x_{1}}\rangle,\ldots,\frac{u}{\vert u\vert}\langle u,u_{x_{n}}\rangle\right)dx\\
&\leq\max\,\{1,\Vert\Phi'_{k}\Vert_{L^{\infty}(\mathbb{R}_{+})}^{\tau}\}\int_{B_{R}\cap\{a_{k}\,<\,\vert u\vert\,<\,b_{k}\}}g(x,\vert u\vert\,\vert u_{x_{1}}\vert,\ldots,\vert u\vert\,\vert u_{x_{n}}\vert)\,dx\\
&\leq\max\,\{1,\Vert\Phi'_{k}\Vert_{L^{\infty}(\mathbb{R}_{+})}^{\tau}\}\int_{B_{R}\cap\{a_{k}\,<\,\vert u\vert\,<\,b_{k}\}}g(x,b_{k}\,\vert u_{x_{1}}\vert,\ldots,b_{k}\,\vert u_{x_{n}}\vert)\,dx\\
&=\max\,\{1,\Vert\Phi'_{k}\Vert_{L^{\infty}(\mathbb{R}_{+})}^{\tau}\}\int_{B_{R}\cap\{a_{k}\,<\,\vert u\vert\,<\,b_{k}\}}f(x,b_{k}\,Du)\,dx\\
&\leq\max\,\{1,\Vert\Phi'_{k}\Vert_{L^{\infty}(\mathbb{R}_{+})}^{\tau}\}\,\,b_{k}^{\tau}\int_{B_{R}}f(x,Du)\,dx
\end{align*}and the last integral is finite, since $u\in W^{1,\mathcal{F}}(\Omega;\mathbb{R}^{m})$.
To prove that $I_{2}$ is bounded, we apply Lemma \ref{lemma1:daipotesi}
(i) and the definition of $\Phi_{k}$, thus obtaining 
\[
I_{2}=\int_{B_{R}\cap\{\vert u\vert\,>\,a_{k}\}}f(x,\eta^{q}\,\Phi_{k}(\vert u\vert)\,Du)\,dx\,\leq\,k^{\tau}\int_{B_{R}}f(x,Du)\,dx<+\infty\,.
\]
We now turn our attention to $I_{3}$. Using Lemma \ref{lemma1:daipotesi}
(i), the definition of $\Phi_{k}$, ($\mathbf{A2}$), (\ref{eq:eta}),
Lemma \ref{lem:g function} (v) and (\ref{growth}), we get\begin{align*}
I_{3}&\leq\,(qk)^{\tau}\int_{B_{R}}f(x,u\,\eta_{x_{1}},\ldots,u\,\eta_{x_{n}})\,dx\\
&=\,(qk)^{\tau}\int_{B_{R}}g(x,\vert u\vert\,\vert\eta_{x_{1}}\vert,\ldots,\vert u\vert\,\vert\eta_{x_{n}}\vert)\,dx\\
&\leq\,(qk)^{\tau}\int_{B_{R}}g\left(x,\frac{2\vert u\vert}{R-\rho},\ldots,\frac{2\vert u\vert}{R-\rho}\right)dx\\
&\leq\left(\frac{2qk}{R-\rho}\right)^{\tau}\int_{B_{R}}\mu(x)\left\{ 1+n^{\frac{q}{2}}|u|^{q}\right\} dx\,,
\end{align*}where we have also used $\frac{2}{R-\rho}>2$. Since $\vert u\vert\in L_{loc}^{\overline{\sigma}^{*}}(\Omega)$
by Proposition \ref{prop:summability} and $qs'<\overline{\sigma}^{*}$
by assumption (see Remark \ref{thm:qs'}), the last integral in the
preceding estimate is finite. Indeed, by Hölder's inequality,
\[
\int_{B_{R}}\mu(x)\vert u\vert^{q}\,dx\,\leq\,\Vert\mu\Vert_{L^{s}(B_{R_{0}})}\,\Vert u\Vert_{L^{qs'}(B_{R_{0}})}^{q}<+\infty\,.
\]
\foreignlanguage{british}{$\hspace*{1em}$Let u}s now consider the
Euler's equation \foreignlanguage{british}{(\ref{p:eulero'})} with
test function $\varphi_{k}$. We obtain\begin{align}\label{eq:I4I5I6}
I_{4}+I_{5}&:=\sum_{j=1}^{n}\sum_{\alpha=1}^{m}\int_{B_{R}}\frac{\partial f}{\partial{\xi_{j}^{\alpha}}}(x,Du)\,u_{x_{j}}^{\alpha}\,\Phi_{k}(|u|)\,\eta^{q}\,dx\nonumber\\
&\,\,\,\,\,\,\,+\sum_{j=1}^{n}\sum_{\alpha,\beta=1}^{m}\int_{B_{R}}\frac{\partial f}{\partial{\xi_{j}^{\alpha}}}(x,Du)\,u^{\alpha}\,\frac{u^{\beta}}{|u|}\,u_{x_{j}}^{\beta}\,\Phi'_{k}(|u|)\,\eta^{q}\,dx\nonumber\\
&\le q\left|\sum_{j=1}^{n}\sum_{\alpha=1}^{m}\int_{B_{R}}\frac{\partial f}{\partial{\xi_{j}^{\alpha}}}(x,Du)\,\Phi_{k}(|u|)\,u^{\alpha}\,\eta^{q-1}\eta_{x_{j}}\,dx\right|=:I_{6}\,.
\end{align}At this stage, we estimate $I_{4}$, $I_{5}$ and $I_{6}$ separately.\foreignlanguage{british}{}\\
\foreignlanguage{british}{}\\
\foreignlanguage{british}{$\hspace*{1em}$}{\textbf{\textsc{Estimate
of $I_{4}$}}}\foreignlanguage{british}{}\\
\foreignlanguage{british}{}\\
\foreignlanguage{british}{$\hspace*{1em}$}As far as $I_{4}$ is concerned,
we use the convexity of $f(x,\cdot)$, (\ref{growth}) and (\ref{eq:eta}).
Thus,\begin{align}\label{eq:I4}
I_{4}&\ge\int_{B_{R}}[f(x,Du)-f(x,0)]\,\Phi_{k}(|u|)\,\eta^{q}\,dx\nonumber\\
&\geq\int_{B_{R}}f(x,Du)\,\Phi_{k}(|u|)\,\eta^{q}\,dx\,-\int_{B_{R}}\mu(x)\,\Phi_{k}(|u|)\,dx\,.
\end{align}\foreignlanguage{british}{$\hspace*{1em}$}{\textbf{\textsc{Estimate
of $I_{5}$}}}\foreignlanguage{british}{}\\
\foreignlanguage{british}{}\\
\foreignlanguage{british}{$\hspace*{1em}$}We claim that $I_{5}\geq0$.
Indeed, by (\ref{eq:struttura}) and the properties (i) and (v) in
Lemma \ref{lem:g function}, we get
\[
\sum_{j=1}^{n}\sum_{\alpha,\beta=1}^{m}\frac{\partial f}{\partial{\xi_{j}^{\alpha}}}(x,Du)\,u^{\alpha}\,u^{\beta}\,u_{x_{j}}^{\beta}\,=\,\sum_{j=1}^{n}\frac{\partial g}{\partial t_{j}}(x,\vert u_{x_{1}}\vert,\ldots,\vert u_{x_{n}}\vert)\,\frac{\left(\sum_{\alpha=1}^{m}u^{\alpha}\,u_{x_{j}}^{\alpha}\right)^{2}}{|u_{x_{j}}|}\geq0\,.
\]
Thus, by the monotonicity of $\Phi_{k}$ we have
\begin{equation}
I_{5}=\int_{B_{R}}\sum_{j=1}^{n}\frac{\partial g}{\partial t_{j}}(x,\vert u_{x_{1}}\vert,\ldots,\vert u_{x_{n}}\vert)\,\frac{\left(\sum_{\alpha=1}^{m}u^{\alpha}\,u_{x_{j}}^{\alpha}\right)^{2}}{|u_{x_{j}}|\,\vert u\vert}\,\Phi'_{k}(|u|)\,\eta^{q}\,dx\,\ge0\,.\label{eq:I5}
\end{equation}
\foreignlanguage{british}{$\hspace*{1em}$}{\textbf{\textsc{Estimate
of $I_{6}$}}}\foreignlanguage{british}{}\\
\foreignlanguage{british}{}\\
\foreignlanguage{british}{$\hspace*{1em}$}Using again (\ref{eq:struttura}),
the properties (i) and (v) in Lemma \ref{lem:g function} and (\ref{eq:eta}),
we have\begin{align}\label{eq:I6est1}
I_{6}&=\,q\left|\sum_{j=1}^{n}\int_{B_{R}}\frac{\partial g}{\partial t_{j}}(x,\vert u_{x_{1}}\vert,\ldots,\vert u_{x_{n}}\vert)\,\frac{\langle u,u_{x_{j}}\rangle}{\vert u_{x_{j}}\vert}\,\Phi_{k}(|u|)\,\eta^{q-1}\,\eta_{x_{j}}\,dx\right|\nonumber\\
&\le\,q\,\sum_{j=1}^{n}\int_{B_{R}}\frac{\partial g}{\partial t_{j}}(x,\vert u_{x_{1}}\vert,\ldots,\vert u_{x_{n}}\vert)\,\vert u\vert\,\Phi_{k}(|u|)\,\eta^{q-1}\,\vert D\eta\vert\,dx\nonumber\\
&\le\,\frac{2q}{R-\rho}\,\sum_{j=1}^{n}\int_{A_{R,j}^{-}\,\cup\,A_{R,j}^{+}}\frac{\partial g}{\partial t_{j}}(x,\vert u_{x_{1}}\vert,\ldots,\vert u_{x_{n}}\vert)\,\vert u\vert\,\Phi_{k}(|u|)\,\eta^{q-1}\,dx\,,
\end{align}where 
\[
A_{R,j}^{-}:=B_{R}\cap\left\{ \eta\ne0,\,|u_{x_{j}}|\le\,\frac{2qL|u|}{\eta(R-\rho)}\right\} 
\]
and 
\[
A_{R,j}^{+}:=B_{R}\cap\left\{ \eta\ne0,\,|u_{x_{j}}|>\,\frac{2qL|u|}{\eta(R-\rho)}\right\} 
\]
with $L>0$ to be chosen later.\foreignlanguage{british}{}\\
For a.e. $x\in A_{R,j}^{-}$ define $H_{j}(x,\cdot):\mathbb{R}_{+}\rightarrow\mathbb{R}_{+}$,
\[
H_{j}(x,\varrho):=\,g(x,\vert u_{x_{1}}(x)\vert,\ldots,\vert u_{x_{j-1}}(x)\vert,\varrho,\vert u_{x_{j+1}}(x)\vert,\ldots,\vert u_{x_{n}}(x)\vert)\,,
\]
of class $C^{1}$ w.r.t. $\varrho$ due to Lemma \ref{lem:g function}
(i). By the properties (i) and (iv) in Lemma \ref{lem:g function}
and by the assumption $x\in A_{R,j}^{-}$, the following inequality
holds: 
\begin{equation}
\frac{\partial g}{\partial t_{j}}(x,\vert u_{x_{1}}\vert,\ldots,\vert u_{x_{n}}\vert)\,\frac{2q|u|}{\eta(R-\rho)}\,\leq\,\frac{1}{L}\,\frac{\partial H_{j}}{\partial\varrho}\left(x,\frac{2qL|u|}{\eta(R-\rho)}\right)\frac{2qL|u|}{\eta(R-\rho)}\,.\label{eq:estimation1}
\end{equation}
Now, let $\omega>1$ and denote by $\mathsf{e}_{1}$ the vector $(1,0,\ldots,0)$
in $\mathbb{R}^{m}$. By the definition of $H_{j}$, (\ref{eq:struttura}),
($\mathbf{A3}$) and Lemma \ref{lem:g function} (v), for a.e. $x\in A_{R,j}^{-}$
and for every $\varrho\in\mathbb{R}_{+}$ we have\begin{align*}
H_{j}(x,\omega\varrho)\,&=\,f\left(x,\omega\,\frac{u_{x_{1}}(x)}{\omega},\ldots,\omega\,\frac{u_{x_{j-1}}(x)}{\omega},\omega\varrho\,\mathsf{e}_{1}^{T},\omega\,\frac{u_{x_{j+1}}(x)}{\omega},\ldots,\omega\,\frac{u_{x_{n}}(x)}{\omega}\right)\\
&\leq\,\omega^{\tau}f\left(x,\frac{u_{x_{1}}(x)}{\omega},\ldots,\frac{u_{x_{j-1}}(x)}{\omega},\varrho\,\mathsf{e}_{1}^{T},\frac{u_{x_{j+1}}(x)}{\omega},\ldots,\frac{u_{x_{n}}(x)}{\omega}\right)\\
&=\,\omega^{\tau}g\left(x,\frac{\vert u_{x_{1}}(x)\vert}{\omega},\ldots,\frac{\vert u_{x_{j-1}}(x)\vert}{\omega},\varrho,\frac{\vert u_{x_{j+1}}(x)\vert}{\omega},\ldots,\frac{\vert u_{x_{n}}(x)\vert}{\omega}\right)\\
&\leq\,\omega^{\tau}g(x,\vert u_{x_{1}}(x)\vert,\ldots,\vert u_{x_{j-1}}(x)\vert,\varrho,\vert u_{x_{j+1}}(x)\vert,\ldots,\vert u_{x_{n}}(x)\vert)\,=\,\omega^{\tau}\,H_{j}(x,\varrho)\,.
\end{align*}Therefore, we can apply Lemma \ref{delta2nuovo} with $h=H_{j}(x,\cdot)$,
thus obtaining
\begin{equation}
\frac{1}{L}\,\frac{\partial H_{j}}{\partial\varrho}\left(x,\frac{2qL|u|}{\eta(R-\rho)}\right)\frac{2qL|u|}{\eta(R-\rho)}\,\leq\,\frac{\tau}{L}\,H_{j}\left(x,\frac{2qL|u|}{\eta(R-\rho)}\right).\label{eq:estimation2}
\end{equation}
Using again the definition of $H_{j}$, Lemma \ref{lem:g function}
(v) and (\ref{eq:struttura}), we get\begin{align}\label{eq:estimationtrois}
H_{j}\left(x,\frac{2qL|u|}{\eta(R-\rho)}\right)&\leq\,g\left(x,\vert u_{x_{1}}\vert,\ldots,\vert u_{x_{j-1}}\vert,\vert u_{x_{j}}\vert+\,\frac{2qL|u|}{\eta(R-\rho)},\vert u_{x_{j+1}}\vert,\ldots,\vert u_{x_{n}}\vert\right)\nonumber\\
&=\,f\left(x,u_{x_{1}},\ldots,u_{x_{j-1}},\left(\vert u_{x_{j}}\vert+\,\frac{2qL|u|}{\eta(R-\rho)}\right)\mathsf{e}_{1}^{T},u_{x_{j+1}},\ldots,u_{x_{n}}\right).
\end{align}Since\begin{align*}
&\left(u_{x_{1}},\ldots,u_{x_{j-1}},\left(\vert u_{x_{j}}\vert+\,\frac{2qL|u|}{\eta(R-\rho)}\right)\mathsf{e}_{1}^{T},u_{x_{j+1}},\ldots,u_{x_{n}}\right)\\
&=\,\frac{1}{2}\,(2u_{x_{1}},\ldots,2u_{x_{j-1}},2\vert u_{x_{j}}\vert\,\mathsf{e}_{1}^{T},2u_{x_{j+1}},\ldots,2u_{x_{n}})\\
&\,\,\,\,\,\,\,+\,\frac{1}{2}\left(\underbrace{0^{T}}_{\in\,\mathbb{R}^{m}},\ldots,\underbrace{0^{T}}_{\in\,\mathbb{R}^{m}},\frac{4qL|u|}{\eta(R-\rho)}\,\mathsf{e}_{1}^{T},\underbrace{0^{T}}_{\in\,\mathbb{R}^{m}},\ldots,\underbrace{0^{T}}_{\in\,\mathbb{R}^{m}}\right)=:\,\frac{1}{2}\,\mathbf{v}+\frac{1}{2}\,\mathbf{w}\,,
\end{align*}by the convexity of $f(x,\cdot)$ we obtain 
\begin{equation}
f\left(x,u_{x_{1}},\ldots,u_{x_{j-1}},\left(\vert u_{x_{j}}\vert+\,\frac{2qL|u|}{\eta(R-\rho)}\right)\mathsf{e}_{1}^{T},u_{x_{j+1}},\ldots,u_{x_{n}}\right)\leq\,\frac{1}{2}\,f(x,\mathbf{v})\,+\,\frac{1}{2}\,f(x,\mathbf{w})\,.\label{eq:estimation4}
\end{equation}
Of course, using (\ref{eq:struttura}) and ($\mathbf{A3}$) we have
\begin{equation}
f(x,\mathbf{v})\,=\,g(x,2\vert u_{x_{1}}\vert,\ldots,2\vert u_{x_{j}}\vert,\ldots,2\vert u_{x_{n}}\vert)\,=\,f(x,2\,Du(x))\,\le\,2^{\tau}f(x,Du(x))\,.\label{eq:estimation5}
\end{equation}
Let us deal with $f(x,\mathbf{w})$. By (\ref{eq:struttura}), Lemma
\ref{lem:g function} (v) and ($\mathbf{A4}$), we deduce\begin{align}\label{eq:estimationsix}
f(x,\mathbf{w})\,&=\,g\left(x,\underbrace{0,\ldots,0}_{j-1},\frac{4qL|u|}{\eta(R-\rho)},\underbrace{0,\ldots,0}_{n-j}\right)\,\leq\,g\left(x,\frac{4qL|u|}{\eta(R-\rho)},\ldots,\frac{4qL|u|}{\eta(R-\rho)}\right)\nonumber\\
&\leq\,\mu(x)\left\{ 1+n^{\frac{q}{2}}\left[\frac{4qL|u|}{\eta(R-\rho)}\right]^{q}\right\}\,\leq\,(4q\sqrt{n})^{q}\,\mu(x)\left\{ 1+\left[\frac{L\vert u\vert}{\eta(R-\rho)}\right]^{q}\right\}.
\end{align}Without loss of generality, we can assume $L\geq1$ so that $\frac{L}{R-\rho}>1$.
Therefore, collecting (\ref{eq:estimation1})$-$\foreignlanguage{british}{\eqref{eq:estimationsix}
and using $0\leq\eta\leq1$, we get\begin{align}\label{eq:I6est2}
&\frac{2q}{R-\rho}\,\sum_{j=1}^{n}\int_{A_{R,j}^{-}}\frac{\partial g}{\partial t_{j}}(x,\vert u_{x_{1}}\vert,\ldots,\vert u_{x_{n}}\vert)\,\vert u\vert\,\Phi_{k}(|u|)\,\eta^{q-1}\,dx\nonumber\\
&\,\,\,\,\,\,\,\leq\,\frac{2^{\tau-1}\,\tau n}{L}\int_{B_{R}}f(x,Du)\,\Phi_{k}(|u|)\,\eta^{q}\,dx\,+\,\frac{(4q\sqrt{n})^{q}\,\tau nL^{q-1}}{2\,(R-\rho)^{q}}\int_{B_{R}}\mu(x)\left\{ 1+\vert u\vert^{q}\right\} \Phi_{k}(|u|)\,dx\,.
\end{align}Let us now deal with $A_{R,j}^{+}$. By the properties (v), (i) and
(iv) in Lemma \ref{lem:g function} we have\begin{align*}
&g(x,2\vert u_{x_{1}}\vert,\ldots,2\vert u_{x_{j}}\vert,\ldots,2\vert u_{x_{n}}\vert)\\
&\,\,\,\,\,\,\,\geq\,g(x,\vert u_{x_{1}}\vert,\ldots,2\vert u_{x_{j}}\vert,\ldots,\vert u_{x_{n}}\vert)\\
&\,\,\,\,\,\,\,\geq\,g(x,\vert u_{x_{1}}\vert,\ldots,\vert u_{x_{j}}\vert,\ldots,\vert u_{x_{n}}\vert)\,+\,\frac{\partial g}{\partial t_{j}}(x,\vert u_{x_{1}}\vert,\ldots,\vert u_{x_{j}}\vert,\ldots,\vert u_{x_{n}}\vert)\,\vert u_{x_{j}}\vert\\
&\,\,\,\,\,\,\,\geq\,\frac{\partial g}{\partial t_{j}}(x,\vert u_{x_{1}}\vert,\ldots,\vert u_{x_{j}}\vert,\ldots,\vert u_{x_{n}}\vert)\,\vert u_{x_{j}}\vert\,.
\end{align*}Using the above inequality together with (\ref{eq:struttura}) and
}($\mathbf{A3}$)\foreignlanguage{british}{, for a.e.} $x\in A_{R,j}^{+}$
\foreignlanguage{british}{we obtain\begin{align*}
\frac{\partial g}{\partial t_{j}}(x,\vert u_{x_{1}}\vert,\ldots,\vert u_{x_{n}}\vert)\,\frac{2q|u|}{\eta(R-\rho)}\,&\leq\,\frac{1}{L}\,\frac{\partial g}{\partial t_{j}}(x,\vert u_{x_{1}}\vert,\ldots,\vert u_{x_{n}}\vert)\,\vert u_{x_{j}}\vert\\
&\leq\,\frac{1}{L}\,g(x,2\vert u_{x_{1}}\vert,\ldots,2\vert u_{x_{j}}\vert,\ldots,2\vert u_{x_{n}}\vert)\\
&=\,\frac{1}{L}\,f(x,2\,Du(x))\,\leq\,\frac{2^{\tau}}{L}\,f(x,Du(x))\,.
\end{align*}Thus\begin{align}\label{eq:I6est3}
\frac{2q}{R-\rho}\,\sum_{j=1}^{n}\int_{A_{R,j}^{+}}\frac{\partial g}{\partial t_{j}}(x,\vert u_{x_{1}}\vert,\ldots,\vert u_{x_{n}}\vert)\,\vert u\vert\,\Phi_{k}(|u|)\,\eta^{q-1}\,dx\,\leq\,\frac{2^{\tau}n}{L}\int_{B_{R}}f(x,Du)\,\Phi_{k}(|u|)\,\eta^{q}\,dx\,,
\end{align}and joining estimates \eqref{eq:I6est1}, \eqref{eq:I6est2} and \eqref{eq:I6est3},
we get\begin{align}\label{eq:I6est4}
I_{6}&\leq\,\frac{2^{\tau}n(\tau+1)}{L}\int_{B_{R}}f(x,Du)\,\Phi_{k}(|u|)\,\eta^{q}\,dx\nonumber\\
&\,\,\,\,\,\,\,+\,\frac{(4q\sqrt{n})^{q}\,\tau nL^{q-1}}{2\,(R-\rho)^{q}}\int_{B_{R}}\mu(x)\left\{ 1+\vert u\vert^{q}\right\} \Phi_{k}(|u|)\,dx\,.
\end{align}Collecting \eqref{eq:I4I5I6}, \eqref{eq:I4}, (\ref{eq:I5}) and
\eqref{eq:I6est4}, and using $\frac{1}{R-\rho}>1$, we find\begin{align}\label{eq:estimationsept}
\int_{B_{R}}f(x,Du)\,\Phi_{k}(|u|)\,\eta^{q}\,dx\,&\leq\,\frac{2^{\tau}n(\tau+1)}{L}\int_{B_{R}}f(x,Du)\,\Phi_{k}(|u|)\,\eta^{q}\,dx\nonumber\\
&\,\,\,\,\,\,\,+\,\frac{(4q\sqrt{n})^{q}\,\tau nL^{q-1}+1}{(R-\rho)^{q}}\int_{B_{R}}\mu(x)\left\{ 1+\vert u\vert^{q}\right\} \Phi_{k}(|u|)\,dx\,.
\end{align}Now we choose $L=2^{\tau+1}\,n(\tau+1)>1$ and reabsorb the first
integral in the right-hand side of \eqref{eq:estimationsept} by the
left-hand side, thus obtaining 
\begin{equation}
\int_{B_{R}}f(x,Du)\,\Phi_{k}(|u|)\,\eta^{q}\,dx\,\le\,\frac{c_{0}}{(R-\rho)^{q}}\int_{B_{R}}\mu(x)\left\{ 1+\vert u\vert^{q}\right\} \Phi_{k}(|u|)\,dx\,,\label{eq:stima1}
\end{equation}
where $c_{0}$ is a positive constant depending only on $n$, $q$
and $\tau$. Inequalities (\ref{growth}) and (\ref{eq:stima1}) imply
\[
\int_{B_{R}}\lambda_{i}(x)\,\vert u_{x_{i}}\vert^{p_{i}}\,\Phi_{k}(\vert u\vert)\,\eta^{q}\,dx\,\leq\,\frac{c_{0}}{(R-\rho)^{q}}\int_{B_{R}}\mu(x)\left\{ 1+\vert u\vert^{q}\right\} \Phi_{k}(|u|)\,dx\,.
\]
At this point, we recall that} $\Phi_{k}=\Phi_{k}^{(i,\gamma)}$\foreignlanguage{british}{.
Using the monotone convergence theorem, we let $k\rightarrow+\infty$
and, by definition of $\Phi_{k}$, we get\begin{align}\label{eq:stima_gamma}
\int_{B_{R}}\lambda_{i}(x)\,\vert u\vert^{p_{i}\gamma}\,\vert u_{x_{i}}\vert^{p_{i}}\,\eta^{q}\,dx\,&\leq\,\frac{c_{0}}{(R-\rho)^{q}}\int_{B_{R}}\mu(x)\left\{ 1+\vert u\vert^{q}\right\} |u|^{p_{i}\gamma}\,dx\nonumber\\
&\leq\,\frac{2c_{0}}{(R-\rho)^{q}}\int_{B_{R}}\mu(x)\left(\max\,\{\vert u\vert,1\}\right)^{q+p_{i}\gamma}dx\,,
\end{align}where, in the last line, we have used
\[
\vert u\vert^{q+p_{i}\gamma}\,+\,\vert u\vert^{p_{i}\gamma}\,\leq\,2\left(\max\,\{\vert u\vert,1\}\right)^{q+p_{i}\gamma}.
\]
$\hspace*{1em}$In a similar way to what has been done up to now,
it is easy to prove that $\varphi:=u\eta^{q}$ is a test function
for the Euler's equation (\ref{p:eulero'}). Moreover, using this
test function and arguing exactly as above, we obtain
\[
\int_{B_{R}}\lambda_{i}(x)\,\vert u_{x_{i}}\vert^{p_{i}}\,\eta^{q}\,dx\,\leq\,\frac{2c_{0}}{(R-\rho)^{q}}\int_{B_{R}}\mu(x)\left(\max\,\{\vert u\vert,1\}\right)^{q}dx\,,
\]
so that \eqref{eq:stima_gamma} also holds for $\gamma=0$.}\\
\foreignlanguage{british}{$\hspace*{1em}$Therefore,} for every $\gamma\geq0$
and every $i\in\{1,\ldots,n\}$ we have\foreignlanguage{british}{
\[
\int_{B_{R}}\lambda_{i}(x)\,\vert u\vert^{p_{i}\gamma}\,\vert u_{x_{i}}\vert^{p_{i}}\,\eta^{q}\,dx\,\leq\,\frac{c}{(R-\rho)^{q}}\int_{B_{R}}\mu(x)\left(\max\,\{\vert u\vert,1\}\right)^{q+p_{i}\gamma}dx\,,
\]
where $c=c(n,q,\tau)>0$. This concludes the proof.\end{proof}}

\selectlanguage{british}%
\noindent $\hspace*{1em}$We are now in a position to prove Theorem
\ref{thm:theo1}.

\noindent \begin{proof}[\bfseries{Proof of Theorem~\ref{thm:theo1}}]Let\foreignlanguage{english}{
$u\in W^{1,\mathcal{F}}(\Omega;\mathbb{R}^{m})$ be a local minimizer
of $(\ref{functional})$} and consider $x_{0}\in\Omega$ and $R_{0}\in(0,1]$
such that $B_{R_{0}}:=B_{R_{0}}(x_{0})\Subset\Omega$. Also fix $0<\rho<R<R_{0}$
and consider a cut-off function $\eta\in C_{c}^{\infty}(B_{R}(x_{0}))$
satisfying (\ref{eq:eta}). To shorten our notation, we now set 
\begin{equation}
G(x):=\,\max\,\{\vert u(x)\vert,1\}\,.\label{eq:G}
\end{equation}
We split the proof into three steps.\foreignlanguage{english}{}\\

\noindent $\hspace*{1em}$\foreignlanguage{english}{\textbf{Step 1.}}
First we prove that, if $\delta\geq1$ and $\vert u\vert^{\delta}\in W^{1,(\sigma_{1},\ldots,\sigma_{n})}(B_{R_{0}})$,
then 
\begin{equation}
\big\Vert G^{\delta}\big\Vert_{L^{\overline{\sigma}^{*}}(B_{\rho})}\leq\,\frac{C\,\delta}{(R-\rho)^{\frac{q}{\overline{\mathbf{p}}}}}\,\Vert\mu\Vert_{L^{s}(B_{R})}^{\frac{1}{\overline{\mathbf{p}}}}\,\Vert G\Vert_{L^{qs'}(B_{R})}^{\frac{q}{\overline{\mathbf{p}}}\,-1}\,\Vert G^{\delta}\Vert_{L^{qs'}(B_{R})}\prod_{i=1}^{n}\Vert\lambda_{i}^{-1}\Vert_{L^{r_{i}}(B_{R})}^{\frac{1}{np_{i}}}\label{eq:step1}
\end{equation}
for a constant $C>1$ depending on $m,n,\mathbf{p},q,\mathbf{r},\tau$
if $\overline{\sigma}<n$, and also on $R_{0}$ if $\overline{\sigma}\geq n$.\\
$\hspace*{1em}$To prove the above inequality, we notice that for
every $\gamma:=\delta-1\geq0$ and every\foreignlanguage{english}{
$i\in\{1,\ldots,n\}$ we have}\begin{align}\label{eq:stima2}
\int_{B_R}\left|\left[(\vert u\vert^{\gamma+1}+1)\,\eta^{q}\right]_{x_{i}}\right|^{p_{i}}\lambda_{i}(x)\,dx\,&\le\,2^{p_{i}-1}\int_{B_R}(q\eta^{q-1}\vert D\eta\vert)^{p_{i}}\,(\vert u\vert^{\gamma+1}+1)^{p_{i}}\,\lambda_{i}(x)\,dx\nonumber\\
&\,\,\,\,\,\,\,+\,2^{p_{i}-1}\,(\gamma+1)^{p_{i}}\int_{B_R}\lambda_{i}(x)\,\eta^{qp_{i}}\,\vert u\vert^{\gamma p_{i}}\,\vert u_{x_{i}}\vert^{p_{i}}\,dx\nonumber\\
&=:J_{1,i}+J_{2,i}\,.
\end{align}To estimate $J_{1,i}$, we observe that 
\[
J_{1,i}\le\,\frac{2^{2q-1}q^{q}}{(R-\rho)^{p_{i}}}\int_{B_{R}}(\vert u\vert^{\gamma+1}+1)^{p_{i}}\,\lambda_{i}(x)\,dx\,\leq\,\frac{2^{3q-1}q^{q}}{(R-\rho)^{q}}\int_{B_{R}}G^{q+\gamma p_{i}}\,\lambda_{i}(x)\,dx\,,
\]
where we have used $p_{i}\leq q$, $\frac{1}{R-\rho}>1$ and
\[
(\vert u\vert^{\gamma+1}+1)^{p_{i}}\,\leq\,(2\,G^{\gamma+1})^{p_{i}}\,\leq\,2^{q}\,G^{q+\gamma p_{i}}\,.
\]
From (\ref{growth}) it easily follows that 
\[
\lambda_{i}\,\leq\,2\,n^{\frac{q}{2}}\,\mu\,\,\,\,\,\,\,\,\,\,\mathrm{a.e.}\,\,\mathrm{in}\,\,\Omega
\]
for every $i\in\{1,\ldots,n\}$. Therefore, 
\begin{equation}
J_{1,i}\le\,\frac{2^{3q}q^{q}n^{\frac{q}{2}}}{(R-\rho)^{q}}\int_{B_{R}}\mu(x)\,G^{q+\gamma p_{i}}\,dx\,.\label{eq:J1i}
\end{equation}
Since $\eta^{qp_{i}}\leq\eta^{q}$, we can estimate $J_{2,i}$ using
inequality (\ref{eq:jolly}). Thus,\begin{align}\label{eq:J2i}
J_{2,i}&\leq\,2^{q-1}\,(\gamma+1)^{p_{i}}\int_{B_{R}}\lambda_{i}(x)\,\eta^{q}\,\vert u\vert^{\gamma p_{i}}\,\vert u_{x_{i}}\vert^{p_{i}}\,dx\nonumber\\
&\leq\,C_{0}\,\frac{(\gamma+1)^{p_{i}}}{(R-\rho)^{q}}\int_{B_{R}}\mu(x)\,G^{q+\gamma p_{i}}\,dx\,,
\end{align}where $C_{0}=C_{0}(n,q,\tau)>0$. By Proposition \ref{prop:Poincar=0000E9}
applied to $v=(\vert u\vert^{\gamma+1}+1)\,\eta^{q}$, we obtain
\[
\left(\int_{B_{R}}[(\vert u\vert^{\gamma+1}+1)\,\eta^{q}]^{\overline{\sigma}^{*}}dx\right)^{\frac{n}{\overline{\sigma}^{*}}}\leq\,C_{1}\,\prod_{i=1}^{n}\Vert\lambda_{i}^{-1}\Vert_{L^{r_{i}}(B_{R})}^{\frac{1}{p_{i}}}\left(\int_{B_{R}}\left|\left[(\vert u\vert^{\gamma+1}+1)\,\eta^{q}\right]_{x_{i}}\right|^{p_{i}}\lambda_{i}(x)\,dx\right)^{\frac{1}{p_{i}}}
\]
for a positive constant $C_{1}$ depending on $m,n,\mathbf{p},\mathbf{r}$
if $\overline{\sigma}<n$, and also on $R_{0}$ if $\overline{\sigma}\geq n$.
Collecting this inequality, \eqref{eq:stima2}, (\ref{eq:J1i}) and
\eqref{eq:J2i}, we get
\[
\left(\int_{B_{R}}[(\vert u\vert^{\gamma+1}+1)\,\eta^{q}]^{\overline{\sigma}^{*}}dx\right)^{\frac{n}{\overline{\sigma}^{*}}}\leq\,C_{2}\,\prod_{i=1}^{n}\Vert\lambda_{i}^{-1}\Vert_{L^{r_{i}}(B_{R})}^{\frac{1}{p_{i}}}\left[\frac{1+(\gamma+1)^{p_{i}}}{(R-\rho)^{q}}\int_{B_{R}}\mu(x)\,G^{q+\gamma p_{i}}\,dx\right]^{\frac{1}{p_{i}}},
\]
where $C_{2}$ is a positive constant depending on $m,n,\mathbf{p},q,\mathbf{r},\tau$
if $\overline{\sigma}<n$, and additionally on $R_{0}$ if $\overline{\sigma}\geq n$.
Using (\ref{eq:eta}) and (\ref{eq:harmonic}), we then find\begin{align}\label{eq:stima4}
&\left(\int_{B_{\rho}}(\vert u\vert^{\gamma +1}+1)^{\overline{\sigma}^{*}}dx\right)^{\frac{n}{\overline{\sigma}^{*}}}\nonumber\\
&\,\,\,\,\,\,\,\leq\,C_{2}\,\prod_{i=1}^{n}\,\frac{\gamma+1}{(R-\rho)^{\frac{q}{p_{i}}}}\,\Vert\lambda_{i}^{-1}\Vert_{L^{r_{i}}(B_{R})}^{\frac{1}{p_{i}}}\left(\int_{B_{R}}\mu(x)\,G^{q+\gamma p_{i}}\,dx\right)^{\frac{1}{p_{i}}}\nonumber\\
&\,\,\,\,\,\,\,\leq\,\frac{C_{2}\,(\gamma+1)^{n}}{(R-\rho)^{\frac{nq}{\overline{\mathbf{p}}}}}\,\prod_{i=1}^{n}\Vert\lambda_{i}^{-1}\Vert_{L^{r_{i}}(B_{R})}^{\frac{1}{p_{i}}}\left(\int_{B_{R}}\mu(x)\,G^{q+\gamma p_{i}}\,dx\right)^{\frac{1}{p_{i}}},
\end{align}for a different constant $C_{2}>0$. By Hölder's inequality and $\mu\in L_{loc}^{s}(\Omega)$,
we obtain
\begin{equation}
\left(\int_{B_{R}}\mu(x)\,G^{q+\gamma p_{i}}\,dx\right)^{\frac{1}{p_{i}}}\leq\,\Vert\mu\Vert_{L^{s}(B_{R})}^{\frac{1}{p_{i}}}\left(\int_{B_{R}}G^{(q+\gamma p_{i})s'}\,dx\right)^{\frac{1}{p_{i}s'}}.\label{eq:stima5}
\end{equation}
$\hspace*{1em}$If $p_{i}=q$ for some $i\in\{1,\ldots,n\}$, the
above estimate gives
\[
\left(\int_{B_{R}}\mu(x)\,G^{(\gamma+1)q}\,dx\right)^{\frac{1}{q}}\leq\,\Vert\mu\Vert_{L^{s}(B_{R})}^{\frac{1}{q}}\,\Vert G^{\gamma+1}\Vert_{L^{qs'}(B_{R})}\,.
\]
$\hspace*{1em}$If $p_{i}<q$ for some $i\in\{1,\ldots,n\}$, we apply
Hölder's inequality again to the last integral in (\ref{eq:stima5}).
We thus obtain
\begin{equation}
\left(\int_{B_{R}}G^{(q+\gamma p_{i})s'}\,dx\right)^{\frac{1}{p_{i}s'}}=\left(\int_{B_{R}}G^{(q-p_{i})s'}\,G^{(\gamma+1)p_{i}s'}\,dx\right)^{\frac{1}{p_{i}s'}}\leq\,\Vert G\Vert_{L^{qs'}(B_{R})}^{\frac{q}{p_{i}}\,-1}\,\Vert G^{\gamma+1}\Vert_{L^{qs'}(B_{R})}\,.\label{eq:stima6}
\end{equation}
Then, up to redefine the constant $C_{2}>0$, estimates \eqref{eq:stima4},
(\ref{eq:stima5}) and (\ref{eq:stima6}) give
\[
\big\Vert\vert u\vert^{\gamma+1}+1\big\Vert_{L^{\overline{\sigma}^{*}}(B_{\rho})}\leq\,\frac{C_{2}\,(\gamma+1)}{(R-\rho)^{\frac{q}{\overline{\mathbf{p}}}}}\,\Vert\mu\Vert_{L^{s}(B_{R})}^{\frac{1}{\overline{\mathbf{p}}}}\,\Vert G\Vert_{L^{qs'}(B_{R})}^{\frac{q}{\overline{\mathbf{p}}}\,-1}\,\Vert G^{\gamma+1}\Vert_{L^{qs'}(B_{R})}\prod_{i=1}^{n}\Vert\lambda_{i}^{-1}\Vert_{L^{r_{i}}(B_{R})}^{\frac{1}{np_{i}}}\,,
\]
where we have used (\ref{eq:harmonic}) again. Since $\delta=\gamma+1$
and $G^{\delta}\leq\vert u\vert^{\delta}+1$, inequality (\ref{eq:step1})
holds true.\\
\\
$\hspace*{1em}$\foreignlanguage{english}{\textbf{Step 2.}} Now we
prove the local boundedness of $G$, and then of $u$, using the Moser
iteration technique. To simplify our notation even more, we introduce
the function 
\begin{equation}
A(r):=\,\Vert\mu\Vert_{L^{s}(B_{r})}^{\frac{1}{\overline{\mathbf{p}}}}\,\Vert G\Vert_{L^{qs'}(B_{r})}^{\frac{q}{\overline{\mathbf{p}}}\,-1}\prod_{i=1}^{n}\Vert\lambda_{i}^{-1}\Vert_{L^{r_{i}}(B_{r})}^{\frac{1}{np_{i}}}\,,\,\,\,\,\,\,\,\,\,\,r\in(0,R_{0}]\,.\label{eq:A}
\end{equation}
Observe that $A(r)$ is non-decreasing, since
\[
1<p_{i}\le q\,\,\,\,\forall\,\,i\in\{1,\ldots,n\}\,\,\,\,\,\,\,\,\,\,\mathrm{implies}\,\,\,\,\,\,\,\,\,\,\frac{q}{\overline{\mathbf{p}}}-1\geq0\,.
\]
Using (\ref{eq:A}), inequality (\ref{eq:step1}) reads as follows:
\begin{equation}
\Vert G\Vert_{L^{\delta\,\overline{\sigma}^{*}}(B_{\rho})}\,\leq\,\left\{ \frac{C\,\delta\,A(R)}{(R-\rho)^{\frac{q}{\overline{\mathbf{p}}}}}\right\} ^{\frac{1}{\delta}}\Vert G\Vert_{L^{\delta qs'}(B_{R})}\,,\label{eq:stima7}
\end{equation}
where $C>1$. For all $h\in\mathbb{N}$, we define\pagebreak{} 
\begin{equation}
\delta_{h}:=\left(\frac{\overline{\sigma}^{*}}{qs'}\right)^{h-1},\,\,\,\,\,\,\,\,R_{h}:=\frac{R}{2}\left(1+\frac{1}{2^{h-1}}\right),\,\,\,\,\,\,\,\,\rho_{h}:=R_{h+1}\,.\label{eq:delta_h}
\end{equation}
Notice that $\delta_{h}$ has been chosen in such a way that $\delta_{1}=1$
and $\delta_{h}\,\overline{\sigma}^{*}=\delta_{h+1}\,qs'$. Moreover,
$\delta_{h}$ diverges to $+\infty$ as $h\rightarrow\infty$, since
$qs'<\overline{\sigma}^{*}$ by assumption. By (\ref{eq:stima7}),
replacing $\delta$, $R$ and $\rho$ with $\delta_{h}$, $R_{h}$
and $\rho_{h}$, respectively, we have that 
\[
G\in L^{\delta_{h}qs'}(B_{R_{h}})\,\,\,\,\,\,\,\,\mathrm{implies}\,\,\,\,\,\,\,\,G\in L^{\delta_{h+1}qs'}(B_{R_{h+1}})\,.
\]
More precisely, the inequality
\[
\Vert G\Vert_{L^{\delta_{h+1}qs'}(B_{R_{h+1}})}\,\leq\,\left\{ \frac{C\,\delta_{h}\,A(R_{h})}{(R_{h}-\rho_{h})^{\frac{q}{\overline{\mathbf{p}}}}}\right\} ^{\frac{1}{\delta_{h}}}\Vert G\Vert_{L^{\delta_{h}qs'}(B_{R_{h}})}
\]
holds true for every $h\in\mathbb{N}$. Since $A(R_{h})\leq A(R)$
for all $h\in\mathbb{N}$, we obtain 
\begin{equation}
\Vert G\Vert_{L^{\delta_{h+1}qs'}(B_{R_{h+1}})}\,\leq\,\left\{ \frac{C\,\delta_{h}\,A(R)}{(R_{h}-\rho_{h})^{\frac{q}{\overline{\mathbf{p}}}}}\right\} ^{\frac{1}{\delta_{h}}}\Vert G\Vert_{L^{\delta_{h}qs'}(B_{R_{h}})}\label{eq:G-iteration}
\end{equation}
for every $h\in\mathbb{N}$. In particular, if $h=1$ we get
\begin{equation}
\Vert G\Vert_{L^{\overline{\sigma}^{*}}(B_{3R/4})}\,\leq\,\frac{4^{\frac{q}{\overline{\mathbf{p}}}}\,C}{R^{\frac{q}{\overline{\mathbf{p}}}}}\,\Vert\mu\Vert_{L^{s}(B_{R})}^{\frac{1}{\overline{\mathbf{p}}}}\,\Vert G\Vert_{L^{qs'}(B_{R})}^{\frac{q}{\overline{\mathbf{p}}}}\prod_{i=1}^{n}\Vert\lambda_{i}^{-1}\Vert_{L^{r_{i}}(B_{R})}^{\frac{1}{np_{i}}}\,.\label{eq:stima8}
\end{equation}
Notice that the right-hand side of (\ref{eq:stima8}) is finite, because
$\vert u\vert\in L^{\overline{\sigma}^{*}}(B_{R_{0}})$ by Proposition
\ref{prop:summability} and $qs'<\overline{\sigma}^{*}$ by assumption.
Now, we define
\begin{equation}
M_{h}:=\left(\int_{B_{R_{h}}}G^{\delta_{h}qs'}\,dx\right)^{\frac{1}{\delta_{h}qs'}},\,\,\,\,\,\,\,\,h\in\mathbb{N},\label{eq:M_h}
\end{equation}
so that inequality (\ref{eq:G-iteration}) turns into
\begin{equation}
M_{h+1}\,\leq\,\left\{ \frac{C\,\delta_{h}\,A(R)\,(2^{\frac{q}{\overline{\mathbf{p}}}})^{h+1}}{R^{\frac{q}{\overline{\mathbf{p}}}}}\right\} ^{\frac{1}{\delta_{h}}}M_{h}\,.\label{eq:stima9}
\end{equation}
Iterating the above estimate, we get
\begin{equation}
M_{h+1}\,\leq\,M_{1}\,\prod_{k=1}^{h}\left\{ \frac{C\,\delta_{k}\,A(R)\,(2^{\frac{q}{\overline{\mathbf{p}}}})^{k+1}}{R^{\frac{q}{\overline{\mathbf{p}}}}}\right\} ^{\frac{1}{\delta_{k}}}\label{eq:stima10}
\end{equation}
for any $h\in\mathbb{N}$. Now we observe that
\[
\sum_{k=1}^{\infty}\frac{1}{\delta_{k}}\,=\,\sum_{k=0}^{\infty}\left(\frac{qs'}{\overline{\sigma}^{*}}\right)^{k}=\,\frac{\overline{\sigma}^{*}}{\overline{\sigma}^{*}-qs'}\,,
\]
so that\pagebreak{} 
\begin{equation}
\prod_{k=1}^{h}\left(\frac{C\,A(R)}{R^{\frac{q}{\overline{\mathbf{p}}}}}\right)^{\frac{1}{\delta_{k}}}=\left(\frac{C\,A(R)}{R^{\frac{q}{\overline{\mathbf{p}}}}}\right)^{\sum_{k=1}^{h}\frac{1}{\delta_{k}}}<\left(\frac{C\,(1+A(R))}{R^{\frac{q}{\overline{\mathbf{p}}}}}\right)^{\sum_{k=1}^{h}\frac{1}{\delta_{k}}}<\left(\frac{C\,(1+A(R))}{R^{\frac{q}{\overline{\mathbf{p}}}}}\right)^{\frac{\overline{\sigma}^{*}}{\overline{\sigma}^{*}-\,qs'}},\label{eq:stima11}
\end{equation}
where we have also used $C>1$ and $R\leq1$. Moreover, by the definition
of $\delta_{h}$ in (\ref{eq:delta_h}), we have 
\begin{equation}
\prod_{k=1}^{h}\delta_{k}^{\frac{1}{\delta_{k}}}=\left(\frac{\overline{\sigma}^{*}}{qs'}\right)^{\sum_{k=1}^{h}\frac{k-1}{\delta_{k}}}<\left(\frac{\overline{\sigma}^{*}}{qs'}\right)^{\sum_{k=1}^{\infty}k\left(\frac{qs'}{\overline{\sigma}^{*}}\right)^{k}}=\left(\frac{\overline{\sigma}^{*}}{qs'}\right)^{\frac{\overline{\sigma}^{*}qs'}{\left(\overline{\sigma}^{*}-\,qs'\right)^{2}}}\label{eq:stima12}
\end{equation}
and
\begin{equation}
\prod_{k=1}^{h}(2^{\frac{q}{\overline{\mathbf{p}}}})^{\frac{k+1}{\delta_{k}}}=\,2^{\frac{q}{\overline{\mathbf{p}}}\left(\frac{\overline{\sigma}^{*}}{qs'}\right)^{2}\sum_{k=2}^{h+1}k\left(\frac{qs'}{\overline{\sigma}^{*}}\right)^{k}}<\,2^{\frac{q}{\overline{\mathbf{p}}}\left(\frac{\overline{\sigma}^{*}}{qs'}\right)^{2}\sum_{k=1}^{\infty}k\left(\frac{qs'}{\overline{\sigma}^{*}}\right)^{k}}=\,2^{\frac{\left(\overline{\sigma}^{*}\right)^{3}}{\overline{\mathbf{p}}s'\left(\overline{\sigma}^{*}-\,qs'\right)^{2}}}.\label{eq:stima13}
\end{equation}
Combining estimates (\ref{eq:stima10})$-$(\ref{eq:stima13}), recalling
the definitions (\ref{eq:A}) and (\ref{eq:M_h}), and letting $h\rightarrow+\infty$,
we obtain
\[
\Vert G\Vert_{L^{\infty}(B_{R/2})}\leq\,c_{1}\left[\frac{1}{R^{\frac{q}{\overline{\mathbf{p}}}}}\left(1+\Vert\mu\Vert_{L^{s}(B_{R})}^{\frac{1}{\overline{\mathbf{p}}}}\,\Vert G\Vert_{L^{qs'}(B_{R})}^{\frac{q}{\overline{\mathbf{p}}}\,-1}\prod_{i=1}^{n}\Vert\lambda_{i}^{-1}\Vert_{L^{r_{i}}(B_{R})}^{\frac{1}{np_{i}}}\right)\right]^{\frac{\overline{\sigma}^{*}}{\overline{\sigma}^{*}-\,qs'}}\left[\int_{B_{R}}G^{qs'}\,dx\right]^{\frac{1}{qs'}},
\]
where $c_{1}>1$ is a constant depending on $m,n,\mathbf{p},q,\mathbf{r},s,\tau$
if $\overline{\sigma}<n$, and also on $R_{0}$ if $\overline{\sigma}\geq n$.
Now, using the fact that $G\geq1$, $\frac{q}{\overline{\mathbf{p}}}-1\geq0$
and $R\leq1$, we deduce
\[
1\leq\left[\fint_{B_{R}}G^{qs'}\,dx\right]^{\frac{1}{qs'}\left(\frac{q}{\overline{\mathbf{p}}}\,-1\right)}=\,\vert B_{R}\vert^{\frac{1}{qs'}\left(1-\,\frac{q}{\overline{\mathbf{p}}}\right)}\,\Vert G\Vert_{L^{qs'}(B_{R})}^{\frac{q}{\overline{\mathbf{p}}}\,-1}\,=\,\frac{\omega_{n}^{\frac{1}{qs'}\left(1-\,\frac{q}{\overline{\mathbf{p}}}\right)}}{R^{\frac{n}{qs'}\left(\frac{q}{\overline{\mathbf{p}}}\,-1\right)}}\,\Vert G\Vert_{L^{qs'}(B_{R})}^{\frac{q}{\overline{\mathbf{p}}}\,-1}
\]
and 
\[
\Vert G\Vert_{L^{qs'}(B_{R})}^{\frac{q}{\overline{\mathbf{p}}}\,-1}\,\leq\,\frac{1}{R^{\frac{n}{qs'}\left(\frac{q}{\overline{\mathbf{p}}}\,-1\right)}}\,\Vert G\Vert_{L^{qs'}(B_{R})}^{\frac{q}{\overline{\mathbf{p}}}\,-1}\,,
\]
where $\omega_{n}$ denotes the $n$-dimensional Lebesgue measure
of the unit ball $B_{1}(0)$. Joining the three previous estimates
and recalling the definition of $G$ in (\ref{eq:G}), we get 
\[
\Vert u\Vert_{L^{\infty}(B_{R/2})}\leq\,\frac{c_{1}}{R^{\vartheta_{3}}}\left[1+\Vert\mu\Vert_{L^{s}(B_{R})}^{\frac{1}{\overline{\mathbf{p}}}}\prod_{i=1}^{n}\Vert\lambda_{i}^{-1}\Vert_{L^{r_{i}}(B_{R})}^{\frac{1}{np_{i}}}\right]^{\vartheta_{1}}\big\Vert\vert u\vert+1\big\Vert_{L^{qs'}(B_{R})}^{\vartheta_{2}}\,,
\]
for a different constant $c_{1}>1$ and
\[
\vartheta_{1}:=\,\frac{\overline{\sigma}^{*}}{\overline{\sigma}^{*}-qs'}\,,
\]
\[
\vartheta_{2}:=\,\vartheta_{1}\left(\frac{q}{\overline{\mathbf{p}}}-1\right)+1=\,\frac{q\,(\overline{\sigma}^{*}-\overline{\mathbf{p}}s')}{\overline{\mathbf{p}}\,(\overline{\sigma}^{*}-qs')}\,,
\]
\[
\vartheta_{3}:=\,\vartheta_{1}\left[\frac{q}{\overline{\mathbf{p}}}+\frac{n}{qs'}\left(\frac{q}{\overline{\mathbf{p}}}-1\right)\right]=\,\frac{\overline{\sigma}^{*}\,[q^{2}s'+n\,(q-\overline{\mathbf{p}})]}{\overline{\mathbf{p}}qs'\,(\overline{\sigma}^{*}-qs')}\,.
\]
\\
We have thus proved that $u\in L^{\infty}(B_{R/2}(x_{0});\mathbb{R}^{m})$
and estimate \eqref{eq:u-bound1}.\\
\\
$\hspace*{1em}$\foreignlanguage{english}{\textbf{Step 3.}} Here we
prove estimate \eqref{eq:u-bound2}. Fix $B_{R}(x_{0})\subset B_{R_{0}}(x_{0})\Subset\Omega$
with $R_{0}\in(0,1]$. Notice that if $Q_{\ell}(x_{0})$ denotes the
cube with edges parallel to the coordinate axes, centered at $x_{0}$
and with side length $2\ell$, then $B_{R/\sqrt{n}}(x_{0})\subseteq Q_{R/\sqrt{n}}(x_{0})\subseteq B_{R}(x_{0})$.\\
$\hspace*{1em}$Define $u_{R}:=\fint_{B_{R}(x_{0})}u\,dx$. Since
$u-u_{R}$ is a local minimizer too, by \eqref{eq:u-bound1} and Hölder's
inequality we have \begin{align}\label{eq:finale1}
&\Vert u-u_{R}\Vert_{L^{\infty}(B_{R/(2\sqrt{n})}(x_{0}))}\nonumber\\
&\,\,\,\,\,\,\,\leq\,\frac{c}{R^{\vartheta_{3}}}\left[1+\Vert\mu\Vert_{L^{s}(B_{R/\sqrt{n}})}^{\frac{1}{\overline{\mathbf{p}}}}\prod_{i=1}^{n}\Vert\lambda_{i}^{-1}\Vert_{L^{r_{i}}(B_{R/\sqrt{n}})}^{\frac{1}{np_{i}}}\right]^{\vartheta_{1}}\left[1+\Vert u-u_{R}\Vert_{L^{\overline{\sigma}^{*}}(B_{R/\sqrt{n}})}\right]^{\vartheta_{2}},
\end{align}where $c=c(m,n,\mathbf{p},q,\mathbf{r},s,\tau,R_{0})>0$. By Proposition
\ref{prop:summability} and Theorem \ref{thm:embedding2}, we get
\[
\Vert u-u_{R}\Vert_{L^{\overline{\sigma}^{*}}(B_{R/\sqrt{n}})}\,\le\,\Vert u-u_{R}\Vert_{L^{\overline{\sigma}^{*}}(Q_{R/\sqrt{n}}(x_{0}))}\,\leq\,c\left[\Vert u-u_{R}\Vert_{L^{1}(B_{R})}+\sum_{i=1}^{n}\Vert u_{x_{i}}\Vert_{L^{\sigma_{i}}(B_{R})}\right].
\]
Now we apply the Poincaré inequality as in \cite[page 185]{cupmarmas0}
(see also \cite[page 84]{cupmarmas1}). Thus we obtain
\[
\Vert u-u_{R}\Vert_{L^{1}(B_{R})}\,\leq\,c\left[1+\sum_{i=1}^{n}\Vert u_{x_{i}}\Vert_{L^{1}(B_{R})}\right],
\]
and combining the two previous estimates, we find
\begin{equation}
\Vert u-u_{R}\Vert_{L^{\overline{\sigma}^{*}}(B_{R/\sqrt{n}})}\,\leq\,c\left[1+\sum_{i=1}^{n}\Vert u_{x_{i}}\Vert_{L^{\sigma_{i}}(B_{R})}\right].\label{eq:finale2}
\end{equation}
At this point, arguing as in the proof of Proposition \ref{prop:summability},
we deduce \begin{align}\label{eq:finale3}
\sum_{i=1}^{n}\Vert u_{x_{i}}\Vert_{L^{\sigma_{i}}(B_{R})}\,&\leq\,\sum_{i=1}^{n}\Vert\lambda_{i}^{-1}\Vert_{L^{r_{i}}(B_{R})}^{\frac{1}{p_{i}}}\left[\int_{B_{R}}f(x,Du)\,dx\right]^{\frac{1}{p_{i}}}\nonumber\\
&\leq\,\left[1+\int_{B_{R}}f(x,Du)\,dx\right]^{\frac{1}{p}}\,\sum_{i=1}^{n}\Vert\lambda_{i}^{-1}\Vert_{L^{r_{i}}(B_{R})}^{\frac{1}{p_{i}}}\,,
\end{align}where $p:=\underset{1\,\leq\,i\,\leq\,n}{\min}\{p_{i}\}$. The final
estimate \eqref{eq:u-bound2} follows by collecting \eqref{eq:finale1},
(\ref{eq:finale2}) and \eqref{eq:finale3}.\end{proof}\vspace{-0.6cm}

\end{singlespace}
\begin{singlespace}

\section{Examples of applicability\label{sec:example}}
\end{singlespace}

\noindent $\hspace*{1em}$\foreignlanguage{english}{In this section,
we provide two examples of integrals to which Theorem \ref{thm:theo1}
is applicable. Before presenting these examples, we need to make some
preliminary considerations.}

\noindent $\hspace*{1em}$\foreignlanguage{english}{Let} \foreignlanguage{english}{$\Omega$
be a bounded open subset of $\mathbb{R}^{n}$, $n\geq2$. For every
$i\in\{1,\ldots,n\},$ let $c_{i}\in(0,+\infty)$, $\kappa_{i}\in[0,+\infty)$,
$p_{i}\in(1,+\infty)$, $\alpha_{i}\in[0,n/p_{i})$ and 
\begin{equation}
\lambda_{i}(x):=\,c_{i}\,\vert x\vert^{\alpha_{i}\,p_{i}},\,\,\,\,\,\,x\in\Omega.\label{eq:lamda_i-example-1}
\end{equation}
For any fixed $i\in\{1,\ldots,n\}$, $\lambda_{i}^{-1}\in L^{r_{i}}(\Omega)$
whenever
\begin{equation}
\begin{cases}
\begin{array}{cc}
r_{i}\in\left[1,\frac{n}{\alpha_{i}\,p_{i}}\right) & \mathrm{if}\,\,\alpha_{i}>0,\\
r_{i}\in[1,\infty]\,\,\,\,\,\,\, & \mathrm{if}\,\,\alpha_{i}=0.
\end{array}\end{cases}\label{eq:r_i}
\end{equation}
Now let $q=\max\,\{p_{i}\}$, $C_{0}=\max\,\{c_{i}\}$, $K_{0}=\max\,\{\kappa_{i}\}$,
$\beta\in[0,n)$ and assume that:\vspace{0.3cm}
}

\selectlanguage{english}%
\noindent $\bullet\,\,\,\,0\in\overline{\Omega}$ if $K_{0}=0$ or
$\beta=0$;\vspace{0.3cm}

\begin{singlespace}
\noindent $\bullet\,\,\,\,0\in\mathbb{R}^{n}\,\backslash\,\Omega$
(possibly $0\in\partial\Omega)$ if $K_{0}>0$ and $\beta>0$.\vspace{0.3cm}

\noindent We now consider the functions $\mu:\Omega\rightarrow[0,\infty)$
and $f:\Omega\times\mathbb{R}{}^{m\times n}\to[0,\infty)$, $m\in\mathbb{N}$,
defined respectively by
\begin{equation}
\mu(x):=\,2^{q-1}\,n\left\{ C_{0}\,\underset{1\,\leq\,j\,\leq n}{\max}\,\vert x\vert^{\alpha_{j}\,p_{j}}\,+\,K_{0}\,\vert x\vert^{-\beta}\right\} \label{eq:mu-example}
\end{equation}
and
\begin{equation}
f(x,\xi):=\sum_{i=1}^{n}\lambda_{i}(x)\,\vert\xi_{i}\vert^{p_{i}}\,+\,\vert x\vert^{-\beta}\sum_{i=1}^{n}\kappa_{i}\,\vert\xi_{i}\vert^{p_{i}}\,.\label{eq:f_example}
\end{equation}
The function $f$ in (\ref{eq:f_example}) trivially satisfies the
assumptions $(\mathbf{A2})$ and $(\mathbf{A3})$ with $\tau=\max\,\{p_{i}\}$.
Moreover, if $p_{i}\geq2$ for every $i\in\{1,\ldots,n\},$ then $f$
also fulfills $(\mathbf{A1})$.\\
\foreignlanguage{british}{$\hspace*{1em}$}We are now ready to give
our two examples, using (\ref{eq:sigma_i}), (\ref{eq:lamda_i-example-1})$-$(\ref{eq:f_example})
and the considerations above.
\end{singlespace}
\begin{example}
\begin{singlespace}
\noindent Let $K_{0}=0$ or $\beta=0$. In this case, the function
$\mu$ in (\ref{eq:mu-example}) belongs to $L^{s}(\Omega)$ for all
$s\in[1,\infty]$. Therefore, if (\ref{eq:r_i}) is in force for every
$i\in\{1,\ldots,n\}$ and if $s=\infty$, the map $f$ in (\ref{eq:f_example})
satisfies $(\mathbf{A4})$ with the weights $\lambda_{i},\mu$ defined
as in (\ref{eq:lamda_i-example-1}) and (\ref{eq:mu-example}). In
light of what has just been said, choosing
\begin{equation}
n=2,\,\,\,\,\,\,p_{1}=2,\,\,\,\,\,\,p_{2}=q=\tau=4,\,\,\,\,\,\,\alpha_{1}=\frac{4}{5}\,,\,\,\,\,\,\,\alpha_{2}=\frac{161}{400}\,,\,\,\,\,\,\,r_{1}=\frac{6}{5}\,,\,\,\,\,\,\,r_{2}=\frac{11}{10}\,,\,\,\,\,\,\,s=\infty\,,\label{eq:parameters1}
\end{equation}
it is easy to verify that the function 
\[
f(x,\xi)\,=\left(c_{1}\,\vert x\vert^{\frac{8}{5}}\,+\,\kappa_{1}\right)\vert\xi_{1}\vert^{2}\,+\left(c_{2}\,\vert x\vert^{\frac{161}{100}}\,+\,\kappa_{2}\right)\vert\xi_{2}\vert^{4},\,\,\,\,\,\,\,\,\,\,x\in\Omega,\,\,\xi\in\mathbb{R}^{m\times2},
\]
and the parameters in (\ref{eq:parameters1}) fulfill all the assumptions
of Theorem \ref{thm:theo1}. Therefore, in the case $n=2$, the main
result of this paper is applicable to integrals of the form 
\[
\mathcal{F}_{1}(v)=\int_{\Omega}\left\{ \left(c_{1}\,\vert x\vert^{\frac{8}{5}}\,+\,\kappa_{1}\right)\vert v_{x_{1}}\vert^{2}\,+\left(c_{2}\,\vert x\vert^{\frac{161}{100}}\,+\,\kappa_{2}\right)\vert v_{x_{2}}\vert^{4}\right\} dx\,.
\]
\end{singlespace}
\end{example}

\begin{example}
\begin{singlespace}
\noindent Let us now consider the case $K_{0}>0$ and $\beta\in(0,n)$.
In this case, the function $\mu$ in (\ref{eq:mu-example}) belongs
to $L^{s}(\Omega)$ for all $s\in(1,\frac{n}{\beta})$. Hence, if
(\ref{eq:r_i}) holds true for every $i\in\{1,\ldots,n\}$ and $1<s<\frac{n}{\beta}$,
the function $f$ in (\ref{eq:f_example}) fulfills $(\mathbf{A4})$
with the weights $\lambda_{i},\mu$ defined by (\ref{eq:lamda_i-example-1})
and (\ref{eq:mu-example}). Thanks to the previous considerations,
choosing
\begin{equation}
n=2,\,\,\,\,p_{1}=2,\,\,\,\,p_{2}=q=\tau=4,\,\,\,\,\alpha_{1}=\frac{4}{5}\,,\,\,\,\,\alpha_{2}=\frac{161}{400}\,,\,\,\,\,r_{1}=\frac{6}{5}\,,\,\,\,\,r_{2}=\frac{11}{10}\,,\,\,\,\,\beta=\frac{1}{5}\,,\,\,\,\,s=5\,,\label{eq:parameters2}
\end{equation}
it is immediate to realize that the function 
\[
f(x,\xi)\,=\left(c_{1}\vert x\vert^{\frac{8}{5}}\,+\,\kappa_{1}\,\vert x\vert^{-\frac{1}{5}}\right)\vert\xi_{1}\vert^{2}\,+\left(c_{2}\,\vert x\vert^{\frac{161}{100}}\,+\,\kappa_{2}\,\vert x\vert^{-\frac{1}{5}}\right)\vert\xi_{2}\vert^{4},\,\,\,\,\,\,\,\,\,\,x\in\Omega,\,\,\xi\in\mathbb{R}^{m\times2},
\]
and the parameters in (\ref{eq:parameters2}) also satisfy all the
assumptions of Theorem \ref{thm:theo1}. Therefore, in the case $n=2$,
our main result can also be applied to integrals of the type
\[
\mathcal{F}_{2}(v)=\int_{\Omega}\left\{ \left(c_{1}\vert x\vert^{\frac{8}{5}}\,+\,\kappa_{1}\,\vert x\vert^{-\frac{1}{5}}\right)\vert v_{x_{1}}\vert^{2}\,+\left(c_{2}\,\vert x\vert^{\frac{161}{100}}\,+\,\kappa_{2}\,\vert x\vert^{-\frac{1}{5}}\right)\vert v_{x_{2}}\vert^{4}\right\} dx\,.
\]
\medskip{}
\end{singlespace}
\end{example}

\selectlanguage{british}%
\begin{singlespace}
\noindent \textbf{Acknowledgements. }The authors are members of the
Gruppo Nazionale per l'Analisi Matematica, la Probabilità e le loro
Applicazioni (GNAMPA) of the Istituto Nazionale di Alta Matematica
(INdAM). P. Ambrosio has been partially supported through the INdAM\textminus GNAMPA
2025 Project ``Regolarità ed esistenza per operatori anisotropi”
(CUP E5324001950001). In addition, P. Ambrosio and G. Cupini acknowledge
financial support under the National Recovery and Resilience Plan
(NRRP), Mission 4, Component 2, Investment 1.1, Call for tender No.
104 published on 2.2.2022 by the Italian Ministry of University and
Research (MUR), funded by the European Union - NextGenerationEU -
Project PRIN\_CITTI 2022 - Title ``Regularity problems in sub-Riemannian
structures'' - CUP J53D23003760006 - Bando 2022 - Prot. 2022F4F2LH.\\

\end{singlespace}

\begin{singlespace}

\lyxaddress{\noindent \textbf{$\quad$}\\
$\hspace*{1em}$\textbf{Pasquale Ambrosio} \\
Dipartimento di Matematica, Università di Bologna\\
Piazza di Porta S. Donato 5, 40126 Bologna, Italy.\\
\textit{E-mail address}: pasquale.ambrosio@unibo.it}

\lyxaddress{\noindent $\hspace*{1em}$\textbf{Giovanni Cupini}\\
Dipartimento di Matematica, Università di Bologna\\
Piazza di Porta S. Donato 5, 40126 Bologna, Italy.\\
\textit{E-mail address}: giovanni.cupini@unibo.it}

\lyxaddress{\noindent $\hspace*{1em}$\textbf{Elvira Mascolo}\\
Dipartimento di Matematica e Informatica ``U. Dini'', Università
degli Studi di Firenze\\
Viale Morgagni 67/A, 50134 Firenze, Italy.\\
\textit{E-mail address}: elvira.mascolo@unifi.it}
\end{singlespace}

\end{document}